\documentclass[a4paper, 10pt]{amsart} 

\usepackage[T1]{fontenc}
\usepackage[utf8]{inputenc}
\usepackage{lmodern}
\usepackage{amsmath,amssymb,amsthm,mathtools}
\usepackage{enumitem}
\usepackage{microtype}
\usepackage[colorlinks=true,linkcolor=blue,citecolor=blue,urlcolor=blue]{hyperref}
\usepackage{tikz-cd} 
\numberwithin{equation}{section}

\newtheorem{thm}{Theorem}[section]
\newtheorem{lemma}[thm]{Lemma}
\newtheorem{prop}[thm]{Proposition}
\newtheorem{cor}[thm]{Corollary}

\newtheorem{introthm}{Theorem}

\theoremstyle{definition}
\newtheorem{rmk}[thm]{Remark}
\newtheorem{defn}[thm]{Definition}
\newtheorem{problem}[thm]{Problem}

\title{\ensuremath{C^{*}}-algebras constructed using acyclic groups} 
\author{Jakub Curda}
\address{Mathematical Institute, Andrew Wiles building, Woodstock Rd, OX2 6GG, Oxford, United Kingdom} 
\email{jakub.curda@maths.ox.ac.uk} 
\author{Shanshan Hua} 
\address{Mathematisches Institut, Fachbereich Mathematik und
Informatik der Universität Münster, Einsteinstrasse 62, 48149 Münster, Germany.} 
\email{shua@uni-muenster.de} 
\author{Raz Slutsky}
\address{Merton College, Merton Street, OX1 4JD, Oxford, United Kingdom} 
\email{raz.slutsky@maths.ox.ac.uk} 
\date{} 

\begin{document}
\begin{abstract}
We show that the free group on countably many generators is not \ensuremath{C^{*}}-superrigid.
Furthermore, we realize each of $\ensuremath{\mathcal{O}}_\infty$ and the Jiang--Su algebra as a crossed product of a commutative \ensuremath{C^{*}}-algebra. Finally, we construct non-trivial, outer actions on classifiable \ensuremath{C^{*}}-algebras whose crossed products exhaust all classifiable \ensuremath{C^{*}}-algebras.
The common theme in the proofs is the use of torsion-free, acyclic groups. 
\end{abstract}

\thanks{The first author was supported by the Engineering and Physical Sciences Research Council [EP/X026647/1]. The second author was funded by the Deutsche Forschungsgemeinschaft (DFG, German Research Foundation) under Germany’s Excellence Strategy EXC 2044-390685587, Mathematics Münster: Dynamics–Geometry–Structure, by the SFB 1442 of the DFG. Part of this work was done during the workshop ``Classification and regularity of inclusions'' (Odense, 2026), funded by IRFD Grant 1054-00094B. For the purpose of open access, the authors have applied a CC-BY license to any author accepted manuscript arising from this submission. } 

\maketitle 

\section{Introduction} 
One of the most natural sources of operator algebras, dating back to the origins of the subject (\cite{mvn36}), comes from groups and their actions. Given an action of a group $G$ on a measurable or topological space, 
one obtains an induced action on the corresponding commutative operator algebra. 
The associated crossed products
are generated by this commutative algebra, together with unitaries implementing the induced action.  
In particular, when $X$ is a single point, these constructions recover the group von Neumann algebra $L(G)$, and the reduced group \ensuremath{C^{*}}-algebra $\ensuremath{C^{*}_r}(G)$ which captures the representations of $G$ weakly contained in the left-regular representation. 
The study of crossed products has proved to be highly influential (\cite{popa07-rigidity,GPS,BKKO}), and three main themes have been studied extensively over the decades: \emph{structure} - 
how properties of the group or the action are reflected in the crossed product; \emph{rigidity} - to what extent the crossed product remembers the underlying group or dynamics; 
and \emph{realization} - which algebras can arise from such constructions. 
In this paper, we solve three problems fitting into these themes. 

A group $G$ is called \emph{\ensuremath{C^{*}}-superrigid} if, whenever $\ensuremath{C^{*}_r}(G)\cong \ensuremath{C^{*}_r}(H)$ for some group $H$, it follows that $G\cong H$. 
This is a \ensuremath{C^{*}}-algebraic analogue of \ensuremath{W^{*}}-superrigidity, defined similarly using group von Neumann algebras. 
The latter notion has been extensively studied and remains a highly active area of research. 
The first examples of i.c.c\text{.} \mbox{\ensuremath{W^{*}}-superrigid} groups were obtained in  \cite{IPV13}, 
and the first i.c.c\text{.} \mbox{\ensuremath{W^{*}}-superrigid} example with property (T) was obtained in \cite{CIOS}. 
Very recently, several breakthrough results have been obtained with the aid of LLMs (\cite{zhou,OpenAInonrigidity,OpenAIrigidity}), showing both strong rigidity and non-rigidity phenomena.
Much earlier, Connes showed (\cite{connes76}) that the group von Neumann algebra of every countable amenable i.c.c\text{.} group is isomorphic to the hyperfinite \ensuremath{\mathrm{II}_{1}} factor $\ensuremath{\mathcal{R}}$. 
Thus, within this class, the group von Neumann algebra retains very little information about the underlying group.
Free groups are not \ensuremath{W^{*}}-superrigid (\cite{dykema93}) and it follows from \cite[Theorem~1.5]{dykemaradulescu} that an infinite free product of non-trivial groups can never be \ensuremath{W^{*}}-superrigid. 
The corresponding question for \ensuremath{C^{*}}-algebras was highlighted in \cite{IPV13}, which pointed out the absence of any known infinite torsion-free group failing \ensuremath{C^{*}}-superrigidity.

Our first main result is the \ensuremath{C^{*}}-algebraic analogue of Dykema--R\u adulescu's result. 
In particular, it provides the first non-trivial, torsion-free examples of groups which are not \ensuremath{C^{*}}-superrigid. 

\begin{introthm}
\label{introthm:general-nonrigidity}
    Suppose $G$ is a free product of infinitely many non-trivial, discrete, countable groups.
    Then $G$ is not \ensuremath{C^{*}}-superrigid. In fact, there are uncountably many pairwise non-isomorphic groups whose reduced group \ensuremath{C^{*}}-algebras are isomorphic to $C_{r}^*(G)$. 
    In particular, the free group $F_{\infty}$ on countably many generators is not \ensuremath{C^{*}}-superrigid. 
\end{introthm} 

We turn to the second problem, concerning dynamical realization of operator algebras. In the von Neumann algebraic setting, the hyperfinite \ensuremath{\mathrm{II}_{1}} factor $\ensuremath{\mathcal{R}}$ plays a central role in Connes' classification of separably acting, injective \ensuremath{\mathrm{II}_{1}} factors: it serves as the tensorial unit, which means that $\ensuremath{\mathcal{M}} \otimes \ensuremath{\mathcal{R}}\cong \ensuremath{\mathcal{M}}$ for every $\ensuremath{\mathcal{M}}$ in this class, and turns out to be the unique such factor (\cite{connes76}). Connes' theorem consequently provides many dynamical realizations of $\ensuremath{\mathcal{R}}$.
Indeed, given any free, ergodic, probability-measure-preserving action of an infinite, discrete, countable, amenable group on a 
standard 
probability space, the associated crossed product is isomorphic to $\ensuremath{\mathcal{R}}$ (\cite[Corollary 6.9]{connes76}). 
In analogy with Connes' work, the Elliott classification programme aims to classify unital, simple, separable and nuclear $C^*$-algebras by \mbox{$K$-theoretical} and tracial invariants. 
The programme has been remarkably successful
(see \cite{KP_classi, Phillips2000, Kirchberg_classi, GLN1, GLN2, TWW} and the survey paper \cite{stuart}), and we call the \ensuremath{C^{*}}-algebras covered by these results  \emph{classifiable}. 
Within this framework, the Jiang--Su algebra $\ensuremath{\mathcal{Z}}$ (\cite{Jiang-Su}) plays a role analogous to that of $\ensuremath{\mathcal{R}}$. The \mbox{$K$-theoretical} and tracial data of $\ensuremath{\mathcal{Z}}$ agree with those of the complex numbers, making tensorial absorption of $\ensuremath{\mathcal{Z}}$ a necessary regularity property in classification.  
Given that $\ensuremath{\mathcal{R}}$ admits many crossed product models, it is therefore natural to ask whether $\ensuremath{\mathcal{Z}}$ admits a concrete crossed product model. 
This is an important special case of \cite[Problem~XLV]{99_open_nuclear}, which asks for dynamical models for all classifiable \ensuremath{C^{*}}-algebras; see Section~\ref{sec:intro-commutative} for further discussion. 

\begin{introthm} 
\label{introthm:Z_realize}  
    There exists a discrete, countable group $G$, a compact, metrizable space $X$, and an 
    action $G\curvearrowright X$ such that 
    \begin{equation} 
        C(X)\rtimes_{r} G\cong \ensuremath{\mathcal{Z}}.
    \end{equation} 
\end{introthm}

A complementary question asks whether every classifiable \ensuremath{C^{*}}-algebra can arise from noncommutative dynamics, in particular, as the crossed product of a classifiable \ensuremath{C^{*}}-algebra by a non-trivial, outer action. 
This question appears as \cite[Problem~XLVI]{99_open_nuclear}, and our third main theorem provides a positive answer. 

\begin{introthm}
\label{introthm:classifiable_realize} 
    Every classifiable \ensuremath{C^{*}}-algebra $A$ can be realized as a crossed product of a non-trivial, outer action of a non-trivial, amenable group on $A$. 
\end{introthm}

The common ingredient underlying these results is the use of acyclic groups. 
A group $G$ is called \emph{acyclic} if its integral group homology vanishes outside of degree zero, or equivalently, if its classifying space \(BG\) has the integral homology of a point.
An early finitely presented example is Higman's group (\cite{Higman1951}), whose acyclicity was observed by Dyer and Vasquez (\cite{DyerVasquez1973}). 
Systematic constructions of acyclic groups were subsequently developed by Baumslag, Dyer, and Heller (\cite{baumslag1980topology}), and by Berrick (\cite{Berrick1989}). 
Acyclic groups have played a key role in resolving the Bass conjecture for amenable groups (\cite{BCM2004}), and are also closely connected with the Kan--Thurston theorem and group-theoretic models of topological constructions (\cite{KanThurston1976,baumslag1980topology}).

For our purposes, the key feature of acyclic groups is the $K$-theory of their reduced group \ensuremath{C^{*}}-algebras (see Lemma~\ref{lem:baum-connes}). 
If a torsion-free group $G$ satisfies the Baum--Connes conjecture - known for a vast class of groups (\cite{HK-BC,oyono2001baum,lafforgue12}) - then the assembly map identifies 
the $K$-theory of the reduced group \ensuremath{C^{*}}-algebra $\ensuremath{C^{*}_r}(G)$ with the topological $K$-homology of the classifying space $BG$ (see \cite{valette} for instance). 
When $G$ is also acyclic, 
the $K$-theory of $\ensuremath{C^{*}_r}(G)$ is the minimal possible: $K_0$ being generated by the class of the unit and $K_1$ vanishing. That is,  
\begin{equation} 
K_{0}(C_{r}^*(G)) \cong \mathbb{Z}\cdot [1]_{0} \quad \text{and} \quad K_{1}(C_{r}^*(G)) = 0. 
\end{equation} 
Hence, through the lens of $K$-theory, $\ensuremath{C^{*}_r}(G)$ is indistinguishable from the complex numbers.
In this way, non-trivial group dynamics can be introduced without altering the $K$-theoretical data. 

\subsection{\ensuremath{C^{*}}-superrigidity} 
\label{sec:intro-superrigidity}

Examples of groups which are not \ensuremath{C^{*}}-superrigid are much sparser than their von Neumann algebraic counterparts.
In the presence of torsion, there are elementary counter-examples - for instance, the two groups of order four have the same \ensuremath{C^{*}}-algebra, and this isomorphism can yield further examples. Therefore, the torsion-free case has attracted most of the attention.
One motivation comes from Kaplansky's
unit conjecture.
This conjecture predicted that every invertible element of the complex group ring $\ensuremath{\mathbb{C}} G$ of a discrete, countable, torsion-free group $G$ is a scalar multiple of a group element, so that $\mathbb{C} G$ uniquely determines $G$, see the introduction of \cite{IPV13}, and \cite{Higman40,GardamComplexUnits}.
Thus \ensuremath{C^{*}}-superrigidity may be regarded as a norm-completed version of the conjecture. 
Despite the conjecture being resolved in the negative by Gardam in \cite{GardamComplexUnits}, the group used there is nevertheless \ensuremath{C^{*}}-superrigid (\cite{ckrtw}). 

The lack of understanding of \ensuremath{C^{*}}-superrigidity of torsion-free groups noted in \cite{IPV13} has sparked new interest in the problem, and many classes of groups have since been shown to be \ensuremath{C^{*}}-superrigid. 
In the amenable setting, the known torsion-free \ensuremath{C^{*}}-superrigid groups include abelian groups (\cite{scheinberg}), and classes closely related to those (\cite{eckhardtraum,omland,ckrtw}).
Non-amenable examples also exist and are almost exclusively obtained through completing to the group von Neumann algebra and proving a \ensuremath{W^{*}}-superrigidity result instead (\cite{ci18,cdad1,cdad2}). 
The only exception is the recent preprint \cite{amcfq}, which constructs groups that are \ensuremath{C^{*}}-superrigid but not
\ensuremath{W^{*}}-superrigid, although the proof still relies heavily on von Neumann
algebraic methods.
In contrast to the varied range of positive results, Theorem~\ref{introthm:general-nonrigidity} provides the first torsion-free groups which are not \ensuremath{C^{*}}-superrigid.

Let us now discuss the strategy behind the proof of Theorem~\ref{introthm:general-nonrigidity} in more detail.
For the purposes of proving non-rigidity, one cannot expect to be able to lift a von Neumann algebraic non-rigidity result to the dense \ensuremath{C^{*}}-subalgebras.
Thus, a \ensuremath{C^{*}}-algebraic argument is required and, in particular, one needs a suitable abstract isomorphism theorem.
For group von Neumann algebras such theorems are available (\cite{connes76,dykemaradulescu}) but classification of \ensuremath{C^{*}}-algebras is typically limited to simple, nuclear algebras.
Reduced \ensuremath{C^{*}}-algebras of non-trivial groups never fit into this framework since for non-trivial groups, simplicity and nuclearity of the reduced \ensuremath{C^{*}}-algebra are incompatible.

The core of the argument is the recent breakthrough of Hirshberg and Phillips (\cite[Theorem~3.1]{hp}), who showed that many infinite free products of \ensuremath{C^{*}}-algebras freely absorb certain inductive limits of 1-dimensional NCCW complexes, and thus provided a new mechanism for constructing isomorphisms between reduced free products, which are typically non-nuclear.
Reduced group \ensuremath{C^{*}}-algebras again do not fit into this framework, as the 1-dimensional NCCW requirement is too restrictive.
The upcoming paper \cite{classification2} will cover the special case of classification of \mbox{*-homomorphisms} in the setting of separable, nuclear domains satisfying the UCT, and simple, monotracial codomains with strict comparison by the unique trace.
This classification theorem, which relies among others on Szab\' o's $KK$-uniqueness theorem (\cite{gabor}), can replace Robert's classification result \cite[Theorem~1.0.1]{robertNCCW} in the proof of \cite[Corollary~3.4]{hp}. 
Consequently, the Hirshberg-Phillips intertwining argument can be generalized so that any unital, separable, nuclear \ensuremath{C^{*}}-algebra satisfying the UCT, admitting a faithful trace $\tau$, and with the $K$-theory of the complex numbers can be freely absorbed: 

\begin{introthm}[A combination of \cite{hp}, \cite{classification} and \cite{gabor}]
    \label{introthm:hp-upgrade}
    Let $\{C_i\}_{i\in \ensuremath{\mathbb{N}}}$, $C_i\neq \ensuremath{\mathbb{C}}$ be  unital, separable \ensuremath{C^{*}}-algebras equipped with faithful traces $\tau_{C_i}\in T(C_i)$ and such that each $C_i$ contains a self-adjoint element with diffuse spectral measure (induced by $\tau_{C_i}$).
    Let $A$ be a unital, separable, nuclear \ensuremath{C^{*}}-algebra satisfying the UCT and equipped with a faithful trace $\tau_A$. 
    Suppose further that $K_0(A)\cong \ensuremath{\mathbb{Z}}\cdot [1]_0$ and $K_1(A)=0$.
    Then
    \begin{equation}
        A*_r \mathop{\scalebox{1.7}{$\;*$}\,\,}\limits_{i\in \ensuremath{\mathbb{N}}} C_i \cong \mathop{\scalebox{1.7}{$\;*$}\,\,}\limits _{i\in\ensuremath{\mathbb{N}}} C_i,
    \end{equation}
    where the reduced free products are taken with respect to the traces $\tau_A$ and $\tau_{C_i}$.
\end{introthm} 

As mentioned previously, the \ensuremath{C^{*}}-algebra $\ensuremath{C^{*}_r}(G)$ of a torsion-free, amenable, acyclic group $G$ has the $K$-theory of the complex numbers and so Theorem~\ref{introthm:hp-upgrade} can be applied with $A=\ensuremath{C^{*}_r}(G)$.
Thus, for any torsion-free, amenable, acyclic group $G$ and any infinite free product group $H$, there is an isomorphism
\begin{equation}
    \ensuremath{C^{*}_r}(G*H)\cong \ensuremath{C^{*}_r}(H).
\end{equation}

Finally, the choice of the torsion-free, amenable, acyclic group is important to ensure that $G*H$ and $H$ are non-isomorphic.
We prove an ``acyclic embedding theorem'' for torsion-free, amenable groups (compare with \cite[Example~3.20]{fournier-facio} which does not apply to all amenable groups and produces uncountable groups): 

\begin{introthm}
\label{introthm:embedings-into-acyclic}
    Every discrete, countable, torsion-free, amenable group embeds into a
    countable, discrete, torsion-free, amenable, acyclic group.
\end{introthm} 

This theorem supplies a sufficiently
rich family of acyclic groups to ensure that the groups \(G*H\) and \(H\) are
 non-isomorphic. 
Furthermore, we show that the failure of \ensuremath{C^{*}}-superrigidity is as bad as it gets: For any infinite free product group, there are uncountably many pairwise non-isomorphic groups with isomorphic \ensuremath{C^{*}}-algebras (Theorem~\ref{thm:general-nonrigidity}).

\subsection{Commutative dynamics for strongly self-absorbing \ensuremath{C^{*}}-algebras}
\label{sec:intro-commutative} 
As recalled above, we call unital, simple, separable, nuclear,
$\ensuremath{\mathcal{Z}}$-stable \ensuremath{C^{*}}-algebras satisfying the UCT \emph{classifiable}.
Such algebras are classified by $K$-theoretical and tracial invariants.  
Kirchberg's dichotomy 
(\cite{BK_dichotomy}) divides 
simple, \ensuremath{\mathcal{Z}}-stable algebras into two fundamentally different subclasses: purely infinite algebras, or \emph{Kirchberg algebras}, and stably finite algebras.

The availability of classification results
motivates the search for concrete models of these algebras, which reveal structural properties and symmetries. 
This is illustrated by Blackadar's dynamical realizations of the CAR algebra and its symmetries (\cite{Blackadar_CAR}), and by Spielberg's graph-groupoid approach to lifting prime-order automorphisms of $K$-theory to automorphisms of Kirchberg algebras (\cite{Spielberg_prime_order}). 
Groupoid models are known in many cases: Spielberg and Clark--Fletcher--an Huef constructed groupoid models for Kirchberg algebras (\cite{Spielberg_graph_model,CFaH_ample_Kirchberg}), Deeley--Putnam--Strung obtained an {\'e}tale groupoid model for $\ensuremath{\mathcal{Z}}$ from point-like minimal systems (\cite{DPS-jiang-su}), and Li's Cartan subalgebras result gives twisted groupoid models for stably finite classifiable algebras via Renault's reconstruction theorem (\cite{Li_Cartan,Renault_Cartan}). A more restrictive question is whether classifiable $C^*$-algebras can be further realized as crossed products, or equivalently, as transformation groupoid $C^*$-algebras. This leads to the following problem. 

\begin{problem}[{\cite[Problem XLV]{99_open_nuclear}}] 
\label{introproblem:crossed_product}  
When is a classifiable $C^*$-algebra of the form ${C(X)\rtimes_{r} G}$ for a suitable action $G \curvearrowright X$ of a  discrete, countable group on a compact, Hausdorff space? 
\end{problem} 

Among classifiable \ensuremath{C^{*}}-algebras, strongly self-absorbing \ensuremath{C^{*}}-algebras form a natural test class for Problem \ref{introproblem:crossed_product}. 
A complete list of strongly self-absorbing $C^*$-algebras satisfying the UCT is given as follows: UHF-algebras of infinite type, the Jiang--Su algebra $\mathcal{Z}$, the Cuntz algebras $\mathcal{O}_{2}$ and $\mathcal{O}_{\infty}$, and tensor products of $\mathcal{O}_{\infty}$ with UHF-algebras of infinite type (\cite[Example 1.14]{SSA}, \cite[Corollary~6.7]{TWW}). 
Several cases are already known: 
UHF algebras of infinite type arise from free, minimal actions of locally finite groups on the Cantor set (\cite[Example 8.1.24]{GKPT}), and $\mathcal{O}_{2}$ is the crossed product of a minimal, topologically free action of $(\mathbb{Z}/2\mathbb{Z})\ast (\mathbb{Z}/3\mathbb{Z})$ on the Cantor set (\cite[Example 4.4.7]{RordamStormer}). 
Moreover, a crossed product model for $\mathcal{O}_{\infty}$ gives models for $\mathcal{O}_{\infty} \otimes M_{\mathfrak{p}}$ by taking the product of the corresponding dynamical systems. Thus the remaining basic cases are $\mathcal{Z}$ and $\mathcal{O}_{\infty}$, resolving which would settle Problem \ref{introproblem:crossed_product} for strongly self-absorbing $C^*$-algebras satisfying the UCT. 

We begin with the Jiang--Su algebra $\mathcal{Z}$, and outline the strategy towards Theorem~\ref{introthm:Z_realize}. Note that $\ensuremath{\mathcal{Z}}$ is the unique classifiable \ensuremath{C^{*}}-algebra with $K_{*}(\mathcal{Z}) = (\ensuremath{\mathbb{Z}}\cdot [1]_0 , 0)$ and a unique trace. Ma and Wu proved that if $\mathcal{Z} \cong C(X)\rtimes_{r} G$, then the invariants of $\mathcal{Z}$ impose strong restrictions on the acting group: $G$ must be torsion-free, amenable, and rationally acyclic (\cite[Theorem C]{mawu}).
Starting with the free, minimal, uniquely ergodic action of $\ensuremath{\mathbb{Z}} $ on a finite-dimensional, point-like,
compact, metrizable space constructed by Deeley--Putnam--Strung (\cite{DPS-jiang-su}), we
inductively obtain a free, minimal, uniquely ergodic action of an acyclic group on the infinite product of this space.
By the going-down principle (\cite[Theorem~3.5.1]{KKG-book}), the Elliott invariant of the resulting crossed product is the same as $\ensuremath{\mathcal{Z}}$; see Section~\ref{sec:model-for-Z}.

The main difficulty is proving $\ensuremath{\mathcal{Z}}$-stability of the resulting crossed
product. The final point-like space in
Theorem~\ref{introthm:Z_realize} is infinite-dimensional, so the
finite-dimensional result of Kerr--Naryshkin (\cite{KerrNaryshkin,glasnerliu}) does not apply directly.
At the first stage, however, the underlying space is finite-dimensional,
and the corresponding crossed product is therefore $\ensuremath{\mathcal{Z}}$-stable by \cite{tomswinter}. 
Using the abundance of commuting subgroups in the construction of our acyclic group, each stage commutes with a copy of $\ensuremath{\mathcal{Z}} $.
Passing to the limit gives a unital copy of $\ensuremath{\mathcal{Z}}$ in the central sequence algebra, and hence $\ensuremath{\mathcal{Z}}$-stability. 

The $\ensuremath{\mathcal{O}}_{\infty}$ case requires a different strategy. The acting group must be nonamenable, since an action of an amenable group on a compact space admits an invariant probability measure, which induces a trace on the crossed product, whereas $\ensuremath{\mathcal{O}}_{\infty}$ is traceless. 
We obtain a particular crossed product model by considering actions of acyclic, hyperbolic groups on the Hilbert cube.
Unlike other realization results (\cite{Spielberg_graph_model, CFaH_ample_Kirchberg, Wu}) where the base space is totally disconnected, the Hilbert cube is contractible.

\begin{introthm} 
\label{introthm:O_infty} 
    For every discrete, countable, torsion-free, non-elementary hyperbolic, acyclic group $H$, there is a minimal, topologically free,  topologically amenable action of $H$ on the Hilbert cube $Q$  such that 
    \begin{equation} 
        C(Q)\rtimes_{r} H\cong \ensuremath{\mathcal{O}}_\infty.
    \end{equation} 
\end{introthm}

We begin with the push-forward action of $H$ on the space of Borel probability measures on its Gromov boundary, which in our setting is homeomorphic to the Hilbert cube $Q$. 
Since $Q$ is contractible and hyperbolic groups satisfy the Baum--Connes conjecture (\cite{lafforgue12, MineyevYu2002}), the going-down principle (\cite[Theorem~3.5.1]{KKG-book}) implies that the $K$-theory of the crossed product is the same as that of $\ensuremath{\mathcal{O}}_{\infty}$. 
The difficulty is minimality, since the
original push-forward action is not minimal. A Baire-category argument in the closure of its conjugacy class produces a minimal action, after which Kirchberg--Phillips' classification leads to the desired isomorphism \cite{KP_classi}. 
Pure infiniteness follows from \cite[Theorem~C]{GardellaGeffenKranzNaryshkin2023}.

A direct consequence of Theorem~\ref{introthm:Z_realize} and Theorem~\ref{introthm:O_infty} is the following: 

\begin{introthm} 
\label{introthm:SSA_crossed_prod} 
    Every strongly self-absorbing \ensuremath{C^{*}}-algebra satisfying the UCT can be realized as a commutative crossed product. 
\end{introthm} 

Lastly, it is worth noting that 
Theorem~\ref{introthm:Z_realize} reduces future crossed product realization problems for classifiable \ensuremath{C^{*}}-algebras to a more focused task: realizing the desired invariant from the dynamics, after which tensoring with the crossed product model of $\ensuremath{\mathcal{Z}}$ given in Theorem~\ref{introthm:Z_realize} provides the required regularity.

\subsection{Non-commutative dynamics} 
Rather than requiring classifiable \ensuremath{C^{*}}-algebras to arise from commutative dynamics, one may allow actions on non-commutative \ensuremath{C^{*}}-algebras. 
Non-commutative dynamical models for stabilizations of classifiable \ensuremath{C^{*}}-algebras were previously known.  
Every stabilization of Kirchberg algebras can be realized as a crossed product of a simple, separable, stable A$\mathbb{T}$-algebra of real rank zero by $\ensuremath{\mathbb{Z}}$ (\cite[Proposition 4.3.3(ii)]{RordamStormer}). 
In the stably finite setting, Elliott, Sato and Thomsen realized the stabilization of the Jiang--Su algebra $\ensuremath{\mathcal{Z}}\otimes \ensuremath{\mathcal{K}}$ to be the crossed product $(E\otimes \ensuremath{\mathcal{K}})\rtimes_{r, \alpha} \ensuremath{\mathbb{Z}}$ of an action $\alpha$, where $E$ is a classifiable \mbox{\ensuremath{C^{*}}-algebra} (\cite[Lemma 3.12]{EST_KMS}). 
Consequently, every stable, $\ensuremath{\mathcal{Z}}$-stable \ensuremath{C^{*}}-algebra $A$ admits a non-commutative model, by considering the action $1\otimes \alpha$ on $A\otimes E\otimes \ensuremath{\mathcal{K}}$.  
Beyond intrinsic structural interest, such crossed product models also provide a useful framework for studying flows and KMS states on classifiable \ensuremath{C^{*}}-algebras as demonstrated in \cite{EST_KMS, RN_rr0}.   

Schafhauser, Tikuisis and White formulate the question as follows: 
\begin{problem}[{\cite[Problem XLVI]{99_open_nuclear}}] 
\label{prob:noncom} 
    When can a classifiable \ensuremath{C^{*}}-algebra be written in the form $A\rtimes_{r} G$ for an outer, amenable action of a non-trivial, discrete, countable group $G$ on a classifiable \ensuremath{C^{*}}-algebra $A$?
\end{problem} 

A natural source of actions on non-commutative \ensuremath{C^{*}}-algebras is the Bernoulli shift. 
For a \ensuremath{C^{*}}-algebra $D$ and a countable discrete group $G$, the \emph{Bernoulli shift action} $\beta\colon G \curvearrowright D^{\otimes G}$ permutes the tensor factors by left translation, that is 
\begin{equation} 
\beta_{h} \left(\bigotimes_{g\in G} d_{g}\right) = \bigotimes_{g\in G} d_{h^{-1}g}, \quad h\in G.  
\end{equation} 
on elementary tensors. 
These actions have been studied extensively from the perspectives of group amenability and $K$-theory computations for crossed products (see \cite{CEKN_Bernoulli}; see also \cite[Example 9.11]{classification} for a detailed treatment of \mbox{$\ensuremath{\mathcal{Z}}^{\otimes\ensuremath{\mathbb{Z}}}\rtimes_{r} \ensuremath{\mathbb{Z}}$}). 

When $D$ is a strongly self-absorbing \ensuremath{C^{*}}-algebra, the isomorphism $D^{\otimes G} \cong D$ allows one to view the Bernoulli shift as an action on $D$ itself. 
It is well-known that every Bernoulli action of a non-trivial, countable group on a tracial, strongly self-absorbing $C^*$-algebra is strongly outer; see the proof of Theorem~\ref{thm:outer-action-on-Z} for instance. 
Bernoulli actions on strongly self-absorbing \ensuremath{C^{*}}-algebras provide model actions in Szab{\'o}'s theory of (semi)-strongly self-absorbing $C^*$-dynamical systems. 
Let $\mathcal{C}$ be the smallest class of groups containing the trivial group and closed under isomorphism, countable directed unions, and extensions by $\ensuremath{\mathbb{Z}}$. 
For any $G\in \mathcal{C}$, the Bernoulli shift action of $G$ on $D$ 
is the unique strongly outer action on $D$ up to (very strong) cocycle conjugacy (\cite[Theorem C]{Gabor_actions_SSA}). 

Our main contribution to Problem~\ref{prob:noncom} is to combine the Bernoulli action with a non-trivial, countable, torsion-free, amenable, acyclic group to obtain a crossed product model for $\ensuremath{\mathcal{Z}}$.

\begin{introthm} 
\label{introthm:bernoulli} 
    Let $G$ be an amenable, countable, torsion-free and acyclic group, and let $\beta$ be the Bernoulli shift action of $G$ on $\ensuremath{\mathcal{Z}}$ through an isomorphism $\ensuremath{\mathcal{Z}}^{\otimes G} \cong \ensuremath{\mathcal{Z}}$. 
    Then 
    \begin{equation}
        \ensuremath{\mathcal{Z}}\rtimes_{r,\beta}G\cong \ensuremath{\mathcal{Z}}. 
    \end{equation} 
\end{introthm} 

Strong outerness of the Bernoulli shift on $\ensuremath{\mathcal{Z}}$, together with amenability of $G$, ensures that $\ensuremath{\mathcal{Z}}\rtimes_{r, \beta} G$ is simple, nuclear, $\ensuremath{\mathcal{Z}}$-stable, and has a unique trace. 
Acyclicity of $G$, combined with the going-down principle, implies that the crossed product has the same $K$-theory as $\ensuremath{\mathcal{Z}}$. 
Classification then identifies the resulting crossed product with $\ensuremath{\mathcal{Z}}$. 
The non-trivial, amenable, countable, torsion-free and acyclic group $\mathfrak{B}(\ensuremath{\mathbb{Z}})$ constructed in this paper belongs to $\mathcal{C}$. 
For this group, Szab{\'o}'s uniqueness theorem therefore identifies the Bernoulli action as the unique strongly outer action on $\ensuremath{\mathcal{Z}}$, up to (very strong) cocycle conjugacy. 

The crossed product model in Theorem~\ref{introthm:bernoulli} extends to every $\ensuremath{\mathcal{Z}}$-stable \ensuremath{C^{*}}-algebra $A$, by a similar argument following Elliott–Sato–Thomsen. 
Indeed, one can consider the tensor-product action $\ensuremath{\mathrm{id}}\otimes\beta$ of $G$ on $A\otimes\ensuremath{\mathcal{Z}}$ along an isomorphism $A\otimes \ensuremath{\mathcal{Z}} \cong A$. 
The resulting crossed product is isomorphic to $A$. 
Since this argument does not require unitality of $A$, it also covers the potentially non-unital classifiable \ensuremath{C^{*}}-algebras considered by Schafhauser, Tikuisis, and White in \cite[Problem XLVI]{99_open_nuclear}. 
We therefore answer Problem~\ref{prob:noncom} affirmatively in full generality.

The paper is organised as follows.
Section~\ref{sec:prelims} contains background material on acyclic groups, $K$-theory, crossed products of \ensuremath{C^{*}}-algebras, and the Jiang--Su algebra. In Section \ref{sec:construction-of-groups}, we construct countable, torsion-free, acyclic, amenable groups and prove Theorem \ref{introthm:embedings-into-acyclic}. Section \ref{section:uncountable_family} proves that given a countable group, there are uncountably many such amenable, acyclic groups such that the free products with the given group remain non-isomorphic. Section \ref{section:classification} contains the needed classification results and a proof of Theorem \ref{introthm:hp-upgrade}. In Section \ref{sec:counterexample} we combine the results obtained thus far to prove Theorem \ref{introthm:general-nonrigidity}. In Section \ref{sec:model-for-Z} we prove the dynamical realisation needed to prove Theorem \ref{introthm:Z_realize}, and in Section \ref{sec:oinfty} we prove Theorem \ref{introthm:O_infty}. Finally, in Section \ref{sec:appli_XLVI} we prove Theorem \ref{introthm:bernoulli} and deduce a solution to Problem \ref{prob:noncom}.

\section*{Developments after 6th October 2026} 
There have been several developments since the OpenAI release on the 6th of October 2026.
The resolution of the free group factor problem (\cite{openai-fgfp}) has led to similar advances in \ensuremath{C^{*}}-algebra theory.
Shlyakhtenko (\cite[Theorem~5.2]{shlyakhtenko-semicircular}) obtained free absorption of $C[0,1]$ into certain finite free products, thus generalizing the Hirshberg--Phillips theorem (\cite[Theorem~3.1]{hp}).
Our techniques, therefore, now extend to obtain further examples of discrete, countable groups which are not \ensuremath{C^{*}}-superrigid.
Indeed, any free product of at least two infinite, discrete, countable groups is not \ensuremath{C^{*}}-superrigid, as we briefly explain next.
In particular, free groups $F_{n}$ with $n\geq 2$ are not
\ensuremath{C^{*}}-superrigid, which answers \cite[Problem~23]{OWR2025Cstar}. 

Without loss of generality, we may only consider a free product $G_1*G_2$ of two infinite groups.
By \cite[Corollary~5.4]{shlyakhtenko-semicircular}, 
\begin{equation}
    \ensuremath{C^{*}_r}(G_1*G_2)\cong \ensuremath{C^{*}_r}(G_1*G_2)*_rC[0,1]^{*\infty} ,
\end{equation}
where the reduced free products are taken with respect to the canonical trace on $\ensuremath{C^{*}_r}(G_1*G_2)$ and the semicircular measure on the interval $[0,1]$.
It follows from the generalized Hirshberg--Phillips theorem (Theorem~\ref{introthm:hp-upgrade}) and Lemma~\ref{lem:baum-connes} that \mbox{$\ensuremath{C^{*}_r}(G_1*G_2)\cong \ensuremath{C^{*}_r}(G_1*G_2*H)$} for any torsion-free, amenable, acyclic group $H$.
Non-isomorphism of $G_1*G_2$ and $G_1*G_2*H$ can be arranged as in the proof of Theorem~\ref{thm:general-nonrigidity}.

\section*{Acknowledgements}
We would like to thank S. White for explaining to us the state-of-the-art classification results for *-homomorphisms and helpful comments on earlier versions of the manuscript.
We would also like to thank J. Carri\' on, J. Gabe, C. Schafhauser, A. Tikuisis and S. White for allowing us to use their unpublished work in the form of Theorem~\ref{thm:classification}, and C. Schafhauser for pointing out a missing assumption in Theorem~\ref{thm:classification} and suggesting useful references.

\section*{AI disclosure}
ChatGPT was used for literature search and for proofreading. It also identified gaps in Lemma~\ref{lem:baum-connes} and Lemma~\ref{lem:nonamenable-acyclic-point-like-perturbation} and suggested corrections for both.
The strategy of proofs, all other mathematical arguments, and all presentation of the work were done by the human authors.

\section{Preliminaries}
\label{sec:prelims}
Throughout, all groups will be assumed discrete and countable.
In this section, we provide some group-theoretical background with a focus on amenable groups, acyclic groups and free products of groups and \ensuremath{C^{*}}-algebras. An overview of $K$-theoretical results used is also given. 

\subsection{Acyclic groups} 
The definition of an acyclic group $G$ is phrased using the homology groups $H_n(G;\ensuremath{\mathbb{Z}})$. 
We refer the reader to the textbook \cite{weibel} for the basics of group homology. 

\begin{defn} 
    A group $G$ is called \emph{acyclic} if 
    \begin{equation}
        H_n(G;\ensuremath{\mathbb{Z}})=0, \qquad n\geq 1.
    \end{equation}
\end{defn} 

Note that the 0-th homology group $H_0(G;\ensuremath{\mathbb{Z}})$ is always equal to $\ensuremath{\mathbb{Z}}$ by definition. Varadarajan introduced in \cite{Varadarajan1985} the notion of \emph{pseudo-mitosis} as a useful criterion for proving acyclicity. 
An equivalent condition was later introduced independently by Berrick \cite{Berrick1989} under the name \emph{binate}; see also \cite[Definition~4.1]{BCM2004} and \cite[Section 3.9]{Berrick1989}.

\begin{defn}
\label{def:pseudomitosis}
Let $G$ be a group and let \(H\leq G\) be a subgroup. We say that $H$ admits a \emph{pseudo-mitosis} in $G$ if there exist 
two homomorphisms 
\begin{equation}
    \psi_0,\psi_1\colon H\rightarrow G 
\end{equation}
and a fixed element \(g\in G\) such that, for all \(h,h'\in H\),
\begin{align} 
h\psi_1(h)&=\psi_0(h), \label{eq:pm1}\\
[h,\psi_1(h')]&=1, \label{eq:pm2}\\
\psi_1(h)&=g^{-1}\psi_0(h)g. \label{eq:pm3}
\end{align} 
The group \(G\) is \emph{pseudo-mitotic} if every finitely generated subgroup of
\(G\) admits a pseudo-mitosis inside \(G\).
\end{defn}  

Note that the homomorphisms $\psi_{0}$ and $\psi_{1}$ are necessarily injective: by \eqref{eq:pm3}, they have the same kernel, and \eqref{eq:pm1} then forces this kernel to be trivial. Thus, broadly speaking, a pseudo-mitosis of a subgroup $H\leq G$ consists of two commuting copies of $H$ inside $G$ whose product is conjugate to one of the original copies. As shown independently by Varadarajan and Berrick, pseudo-mitotic groups are acyclic. 

\begin{thm}[{\cite[Theorem~1.7]{Varadarajan1985}, \cite[Theorem~1.8]{Berrick1989}}]
\label{thm:pseudomitotic-acyclic}
Every pseudo-mitotic (or equivalently binate) group is acyclic. 
\end{thm}

\subsection{Free products of groups and \ensuremath{C^{*}}-algebras}

This subsection reviews standard results on free products of groups and reduced free products of tracial \ensuremath{C^{*}}-algebras.
We begin with the following fact from combinatorial group theory; we refer the reader to \cite[Section~IV.1]{lyndon-schupp} for further details about free products of groups and reduced words.

\begin{lemma}[{\cite[Theorem~IV.1.6]{lyndon-schupp}}]
    \label{lemma:torsion-free-products}
    The free product $*_{i\in I}G_i$ of torsion-free groups $G_i$ is torsion-free.
\end{lemma}

The second group theoretic ingredient is the following consequence of Kurosh's theorem which identifies maximal amenable subgroups of free products.

\begin{lemma} 
\label{lem:amenable-subgroups-free-product}
Let \(G_1\) and \(G_2\) be groups.

\begin{enumerate}[label=(\roman*)]
    \item Every non-cyclic, torsion-free, amenable subgroup of \(G_1*G_2\) is
    contained in a conjugate of \(G_1\) or of \(G_2\).
    \item If \(G_2\) is non-cyclic, torsion-free and amenable, then \(G_2\) is a
    maximal amenable subgroup of \(G_1*G_2\).
\end{enumerate}
\end{lemma}

\begin{proof}
By the Kurosh subgroup theorem (\cite[Theorem~IV.1.10]{lyndon-schupp}), every subgroup $H$ of \(G_1*G_2\) is a free product of the form
\begin{equation}
    H= F*\left(  \mathop{\scalebox{1.7}{$\;*$}\,\,}\limits_{i\in I} H_{i}\right),
\end{equation}
where $F$ is a free group and each $H_{i}$ is the intersection of $H$ with some conjugate of either $G_1$ or $G_2$. 
A free product of two non-trivial groups contains a non-abelian free subgroup, except when both factors have order two.  
The former case cannot happen if the subgroup $H$ is amenable, and the exceptional case cannot occur if the subgroup $H$ is torsion-free. 
Thus an amenable torsion-free subgroup $H$ is either the rank-one free subgroup $\mathbb{Z}$ or is contained in a conjugate of one factor, proving (i).

For (ii), let \(G_2\leq H\leq G_1*G_2\), with \(H\) amenable.  In the Kurosh
decomposition of \(H\), the subgroup \(H\cap G_2=G_2\) occurs as a free factor.
Thus \(H\cong G_2*H_0\) for some group \(H_0\).  If \(H_0\) were non-trivial,
then \(H\) would contain a non-abelian free subgroup, since
\(G_2\not\cong\ensuremath{\mathbb{Z}}/2\).  Hence \(H_0\) is trivial and \(H=G_2\).
\end{proof}

We now turn to free products of \ensuremath{C^{*}}-algebras.
We only provide a brief overview and refer the reader to \cite[Section~4.7]{brownozawa} for a detailed introduction of reduced free products.
Let $A$ and $B$ be unital \ensuremath{C^{*}}-algebras with distinguished faithful traces $\tau_A$ and $\tau_B$.
The \emph{algebraic free product} $A*_\mathrm{alg}B$ of $A$ and $B$ is the *-algebra spanned by formal words with letters being elements of $A$ and $B$ with multiplication given by word concatenation and the conjugation of a word implemented by conjugating each letter and reversing the order of the word.
The only further relations imposed are the relations from $A$ and $B$ and the identification of the units of $A$ and $B$.

The traces $\tau_A$ and $\tau_B$ give rise to a faithful representation of $A*_\mathrm{alg}B$ on a Hilbert space and the \emph{reduced free product} $(A,\tau_A)*_r (B,\tau_B)$ of $A$ and $B$ is the norm completion of $A*_\mathrm{alg}B$ in this representation.
We warn the reader that the construction depends on the chosen traces.

Given a family of unital \ensuremath{C^{*}}-algebras $\left\{ A_i \right\}_{i\in I} $ with distinguished faithful traces $\tau_i$, one can construct the reduced free product
\begin{equation}
    \mathop{\scalebox{1.7}{$\;*$}\,\,}\limits_{i\in I}(A_i,\tau_i).
\end{equation}

The construction of the reduced free product generalizes the free product of groups in the following sense:
\begin{lemma}[{\cite[Example~4.7.4]{brownozawa}}]
\label{lemma:free-product-correspondnce}
    Let $\left\{ G_i \right\} $ be groups.
    Then, with respect to the canonical traces $\tau_{i}$, there is an isomorphism
    \begin{equation}
        \ensuremath{C^{*}_r}\left(\mathop{\scalebox{1.7}{$\;*$}\,\,}\limits
         _{i\in I}G_i\right)\cong \mathop{\scalebox{1.7}{$\;*$}\,\,}\limits
          _{i\in I} (\ensuremath{C^{*}_r}(G_i),\tau_i).
    \end{equation}
\end{lemma} 
The reduced free product is typically denoted by $*_r$ and the symbol $*$ is reserved for the full free product.
Since no full free products will appear in this paper, we will use the latter symbol for the infinite reduced free product.
Moreover, when the faithful traces are clear from context, they will be omitted from the notation. 

\subsection{Crossed products of \ensuremath{C^{*}}-algebras}
\label{sect:crossed_prod} 
In this section we recall some standard results on crossed products.
See the textbook \cite{williams} for basic definitions.

We begin with standard properties of group actions on topological spaces.
Let $G$ be a discrete, countable group with an action on a compact, Hausdorff space $X$.
The action $G\curvearrowright X$ is called 
\begin{enumerate}[label=(\roman*)] 
\item \emph{minimal}, if every orbit is dense in $X$; 
\item \emph{free}, if for every $g \in G\setminus\{1_{G}\}$ and $x\in X$, we have $g x \neq x$; 
\item \emph{topologically free}, if for every $g \in G\setminus\{1_{G}\}$, the set $\{x\in X \ | \ gx = x\}$ has empty interior; 
\item \emph{uniquely ergodic}, if there is a unique $G$-invariant probability measure on $X$; 
\item \emph{topologically amenable}, if for any $\varepsilon > 0$, there exists a continuous map \mbox{$m\colon X\rightarrow \mathrm{Prob}_{f}(G)$},\footnote{We write $\mathrm{Prob}_{f}(G)$ for the finitely supported probability measures on $G$, equipped with the $\ell^{1}$-norm. The left translation action of $G$ on $\mathrm{Prob}_{f}(G)$ is given by $(g\cdot m) (h) \coloneqq m(g^{-1}\cdot h)$ for $g,h \in G$ and $m\in \mathrm{Prob}_{f}(G)$. We use this formulation of topologically amenable actions in our paper; see \cite[Proposition 2.2]{Claire_amenable} for more details.} with common finite support, such that 
\begin{equation} 
\sup_{x\in X}\| m (gx) - g\cdot m(x) \|_{1} < \varepsilon, \qquad g \in G. 
\end{equation} 
\end{enumerate} 
If the group $G$ is amenable, every action of $G$ is automatically topologically amenable. 

Given a group action, one can form the reduced crossed product $C(X)\rtimes_{r} G$, see \cite[Definition~7.7]{williams} for instance.
The properties of the action are reflected in the structure of the crossed product; we collect relevant results below: 

\begin{prop} 
\label{action_prop} 
Let $G$ be a discrete, countable group and $X$ be a compact, Hausdorff space. Suppose $G$ acts on $X$. Then the following are true for the  associated crossed product $C(X)\rtimes_{r} G$: 
\begin{enumerate}[label=(\roman*)]
\item $C(X)\rtimes_{r} G$ is always unital, and is separable if and only if $X$ is metrizable; 
\item \label{simple}
If the action of $G$ on $X$ is topologically free and minimal, the associated reduced crossed product is simple \textup{(}\cite[Corollary on p.~122]{archboldspielberg}\textup{)}; 
\item $C(X)\rtimes_{r} G$ is nuclear if and only if the action of $G$ is topologically amenable \textup{(}\cite[Theorem~3.3.7 and Corollary~6.2.14]{AnantharamanDelarocheRenault2000}\textup{)}. In this case, it has a tracial state precisely when $G$ itself is amenable
\textup{(}\cite[Lemma~2.2 and Remark~2.3]{GardellaGeffenKranzNaryshkin2023}\textup{)}; 
\item \label{unique_trace}
If $G$ is amenable, then $C(X)\rtimes_{r} G$ has a unique trace whenever the action is free and uniquely ergodic \textup{(}\cite[Corollary~9.6]{zeller-meier}\textup{)}. 
\end{enumerate} 
\end{prop}

\subsection{$K$-theory} 
In this section, we collect several standard facts on the operator-algebraic $K$-theory of reduced group \ensuremath{C^{*}}-algebras and crossed products arising from acyclic groups. 
The notation $K_*(A)=(\ensuremath{\mathbb{Z}}\cdot[1]_0,0)$ means that $K_1(A)$ is trivial and $K_0(A)\cong\ensuremath{\mathbb{Z}} $ is generated by the class of the unit.

The vanishing of homology for a group $G$ has consequences for the $K$-theory of its reduced group \ensuremath{C^{*}}-algebra $\ensuremath{C^{*}_r}(G)$ whenever $G$ satisfies the Baum--Connes conjecture. 
If $G$ is additionally torsion-free and acyclic, the inclusion of the unit $\mathbb{C} \xhookrightarrow{} \ensuremath{C^{*}_r}(G)$ induces an isomorphism in $K$-theory, as shown in the next lemma. 
The proof is standard but we include it for completeness. 

\begin{lemma}
\label{lem:baum-connes} 
    Let $G$ be a torsion-free, acyclic group satisfying the Baum--Connes conjecture. The reduced \ensuremath{C^{*}}-algebra of $G$ satisfies 
    \begin{align}
            K_0(\ensuremath{C^{*}_r}(G))= \ensuremath{\mathbb{Z}}\cdot [1]_0 \quad \text{and} \quad 
            K_1(\ensuremath{C^{*}_r}(G))= 0. 
    \end{align} 
\end{lemma} 

\begin{proof} 
Since $G$ satisfies the Baum-Connes conjecture and is furthermore torsion-free, the assembly map gives an isomorphism (\cite[Section~2.1]{valette}) 
\begin{equation}
\label{eq:bc}
    K_*(\ensuremath{C^{*}_r}(G))\cong K_*(BG),
\end{equation} 
where the right hand side is the topological $K$-homology of the classifying space $BG$. It therefore suffices to identify $K_*(BG)$.

Acyclicity of $G$, together with the identification $H_*(G;\ensuremath{\mathbb{Z}})\cong H_*(BG;\ensuremath{\mathbb{Z}})$, implies that the reduced homology satisfies $\widetilde H_*(BG;\ensuremath{\mathbb{Z}})=0$; see
\cite[Theorem~6.10.5]{weibel}. 
Since $BG$ is a path-connected CW-complex, it follows that $\Sigma BG$ is a simply connected CW-complex. 
The suspension isomorphism $\widetilde H_{*+1}(\Sigma BG;\mathbb Z) 
\cong \widetilde H_*(BG;\mathbb Z)$ implies that $\Sigma BG$ has vanishing reduced homology. 
Since $\Sigma BG$ is a simply connected CW-complex, as a consequence of the relative Hurewicz theorem and Whitehead’s theorem (\cite[Corollary 4.33]{Hatcher}), we have that $\Sigma BG$ is contractible. 
By the suspension isomorphism in reduced topological $K$-homology, it follows that
$\widetilde K_*(BG) \cong \widetilde K_{*+1}(\Sigma BG) =0$. Consequently, $K_0(BG)\cong\ensuremath{\mathbb{Z}}$,  and $K_1(BG)=0$. 

Let $\tau$ denote the canonical trace on $\ensuremath{C^{*}_r}(G)$ and let $\tau_*$ denote the induced map,
\begin{equation}
    \tau_*\colon K_0(\ensuremath{C^{*}_r}(G))\rightarrow \ensuremath{\mathbb{R}}.
\end{equation}
Since $G$ satisfies the Baum-Connes conjecture, the image of $\tau_*$ is equal to the subgroup of integers in $\ensuremath{\mathbb{R}}$  (see \cite[Proposition~6.3.1]{valette}).
Since $\tau_*([1]_0)=1$, the unit of $\ensuremath{C^{*}_r}(G)$ must generate $K_0(\ensuremath{C^{*}_r}(G))$.     
\end{proof}

The next theorem follows directly from \cite[Theorem~3.5.1]{KKG-book} and allows us to compute the $K$-theory of crossed products. 

\begin{thm}[{\cite[Theorem~3.5.1]{KKG-book}}] 
\label{thm:KKG} 
    Let $G$ be a torsion-free  group satisfying the Baum-Connes conjecture with coefficients, and let $G$ act on a unital, separable \ensuremath{C^{*}}-algebra $A$ with $K_{0}(A) = \ensuremath{\mathbb{Z}}\cdot[1]_0$ and $K_{1}(A) = 0$. 
    The $K$-theory of the reduced crossed product satisfies 
    \begin{equation} 
      K_*(A\rtimes_{r} G)\cong  K_*(\ensuremath{C^{*}_r}(G))
    \end{equation} 
    and the isomorphism is implemented by the inclusion of $\ensuremath{C^{*}_r}(G)$ into $A\rtimes_{r} G$ and thus preserves the class of the unit. 
\end{thm}

\begin{proof} 
Let $\iota \colon \mathbb{C} \rightarrow A$ be the unital inclusion. As $\iota$ is $G$-equivariant, it induces a standard $G$-equivariant $KK$-$\ensuremath{\mathbb{C}}$-$A$-cycle (see \cite[Definition~2.25]{EKQR}). Let $x = [\iota] \in KK^G(\ensuremath{\mathbb{C}},A)$. Since $G$ is torsion-free, the only compact subgroup of $G$ is the trivial group. Moreover, the map 
\begin{equation} 
\cdot\otimes\operatorname{res} ^G_1(x)\colon K_*(\ensuremath{\mathbb{C}})\rightarrow K_*(A) 
\end{equation} 
is an isomorphism by the hypothesis on $K_*(A)$. 
The going-down principle (\cite[Theorem~3.5.1]{KKG-book}) therefore gives an isomorphism 
\begin{equation} 
\label{eq:J_G} 
\cdot\otimes_{C^{*}_{r}(G)} J_G(x)\colon K_*(\ensuremath{C^{*}_r}(G))\rightarrow K_*(A\rtimes_{r} G). 
\end{equation} 
By \cite[Proposition~3.6]{EKQR}, $J_G(x)$ is the $KK$-class 
induced by the unital inclusion
\begin{equation}
        \iota\rtimes_{r} G\colon \ensuremath{C^{*}_r}(G)\rightarrow A\rtimes_{r} G.
    \end{equation}
Hence the isomorphism (\ref{eq:J_G}) is induced by $\iota\rtimes_{r}G$ and thus preserves the class of the unit.   
\end{proof} 

Combining the theorem above with Lemma~\ref{lem:baum-connes} gives the following. 

\begin{cor} 
\label{cor:K-theory} 
    Let $G$ be a torsion-free, amenable, acyclic group acting on a compact, metrizable space $X$ with $K_{0}(C(X)) = \ensuremath{\mathbb{Z}}\cdot[1]_0$ and $K_{1}(C(X)) = 0$. Then 
    \begin{equation} 
        K_0(C(X)\rtimes_{r} G)=\ensuremath{\mathbb{Z}}\cdot[1]_0 \quad \text{and} \quad  K_1(C(X)\rtimes_{r} G)=0.
    \end{equation} 
\end{cor} 

\subsection{The Jiang--Su algebra and $\ensuremath{\mathcal{Z}}$-stability} 
The Jiang--Su algebra $\ensuremath{\mathcal{Z}}$ was initially constructed in \cite{Jiang-Su}, and has since become a central object in the Elliott classification programme as tensorial absorption of $\ensuremath{\mathcal{Z}}$ is necessary to obtain a complete classification theorem. Recall that a \ensuremath{C^{*}}-algebra $A$ is \emph{$\ensuremath{\mathcal{Z}}$-stable} if $A\otimes \ensuremath{\mathcal{Z}}\cong A$.  

In many examples, $\ensuremath{\mathcal{Z}}$-stability is difficult to verify. The following equivalent characterization of $\ensuremath{\mathcal{Z}}$-stability for unital, separable $C^*$-algebras is useful in practice. Fix a free ultrafilter $\omega$ on $\mathbb{N}$. Denote by $\ell^{\infty}(A)$ the \ensuremath{C^{*}}-algebra of norm bounded sequences in $A$. 
The \emph{norm ultrapower} $A_{\omega}$ of a \ensuremath{C^{*}}-algebra $A$ is defined by 
\begin{equation} 
    A_{\omega} = \ell^{\infty}(A)/ \{(a_n)_{n} \in \ell^{\infty}(A) \ | \lim_{n\to \omega} \|a_n\| = 0\}. 
\end{equation} 
We write elements in $A_\omega$ using representing sequences $(a_n)_n$ from $\ell^\infty(A)$. 
The \emph{central sequence algebra} of $A$ is defined to be the relative commutant $A_{\omega}\cap A'$, where $A$ embeds diagonally in $A_{\omega}$. The following is \cite[Theorem~2.2]{SSA}.

\begin{prop} 
\label{char_Z_stable} 
Let $A$ be a unital and separable \ensuremath{C^{*}}-algebra. Then $A$ is $\ensuremath{\mathcal{Z}}$-stable if and only if there exists a unital embedding $\ensuremath{\mathcal{Z}} \rightarrow A_{\omega} \cap A'$. 
\end{prop} 

The corollary below is an easy observation that will allow us to obtain $\ensuremath{\mathcal{Z}}$-stability. 

\begin{cor} 
\label{Z_commute_lem} 
Let $A = \varinjlim (A_{n}, \varphi_{n})$ be an inductive limit of an
inductive system of unital, separable \ensuremath{C^{*}}-algebras $A_{n}$ and unital connecting maps $\varphi_{n}$. Suppose that for every $n\geq 0$, there is a unital $^*$-homomorphism $\iota_{n}\colon \ensuremath{\mathcal{Z}} \rightarrow A_{n+1}$ such that $\iota_{n}(\ensuremath{\mathcal{Z}})$ commutes with $\varphi_{n}(A_{n})$ in $A_{n+1}$. Then $A$ is $\ensuremath{\mathcal{Z}}$-stable. 
\end{cor} 

\begin{proof} 
Let $\varphi^{(n)} \colon A_{n}\rightarrow A$ be the induced maps. Take the unital embedding 
\begin{equation} 
\iota\colon \ensuremath{\mathcal{Z}} \rightarrow A_{\omega}, \quad d \mapsto (\varphi^{(n+1)}\circ \iota_{n}(d))_{n}. 
\end{equation} 
We claim that $\iota(\ensuremath{\mathcal{Z}})$ commutes with $A$ in $A_{\omega}$. Since $A$ is the dense union of $\varphi^{(n)}(A_{n})$, for every $a\in A$, there exists $a_{n}\in \varphi^{(n)}(A_{n})$ such that $a$ and $(a_{n})_{n}$ represent the same element in $A_{\omega}$. Since in the $n$-th coordinate, $\iota (\ensuremath{\mathcal{Z}})$ commutes with $\varphi^{(n+1)} \circ \varphi_{n} (A_{n}) = \varphi^{(n)} (A_{n})$, it follows that $\iota(\ensuremath{\mathcal{Z}})$ commutes with $(a_{n})_{n}$ in $A_{\omega}$. Therefore $\iota (\ensuremath{\mathcal{Z}})$ commutes with every $a\in A$. Hence, $\iota$ is a unital embedding of $\ensuremath{\mathcal{Z}}$ into $A_{\omega}\cap A'$, which implies that $A$ is $\ensuremath{\mathcal{Z}}$-stable by Proposition~\ref{char_Z_stable}. 
\end{proof}

\section{Construction of torsion-free amenable acyclic groups}
\label{sec:construction-of-groups}

This section contains a construction of torsion-free, amenable, acyclic groups. 
Theorem~\ref{introthm:embedings-into-acyclic} is also proved: 
starting from a countable, torsion-free, amenable group \(G\), we construct a group
\(\ensuremath{\mathfrak{B}}(G)\) containing \(G\), which is again countable, torsion-free, amenable, and
additionally acyclic. This is inspired by \cite{baumslag1980topology}, where, using HNN extensions, it is shown that every countable group embeds into an acyclic group. However, these acyclic groups are far from being amenable. We note that our construction is very similar to the one in \cite[Example 3.20]{fournier-facio}, which, however, does not apply to all amenable groups, and for which the resulting group is uncountable. 

To construct $\ensuremath{\mathfrak{B}}(G)$, we define a new group \(\mathfrak{C}(G)\) which contains a concrete copy of \(G\), a second commuting ``tail copy'', and a shift element that conjugates the tail copy into the product of the two copies.  
This is exactly the mechanism needed for a pseudo-mitosis. Iterating this construction, we obtain the desired acyclic groups.

\subsection{Construction of $\ensuremath{\mathfrak{B}}(G)$} 
Let $G$ be a group. Define 
\begin{equation} 
\mathfrak{D}(G)\coloneqq\left\{x\colon\ensuremath{\mathbb{Z}}\to G \ \Big|
\begin{array}{l} 
x_k \text{ is constant for all sufficiently positive }k \\
x_k \text{ is trivial for all sufficiently negative }k
\end{array} \right\}. 
\end{equation} 

The group operation is pointwise multiplication. 
Note that $\mathfrak{D}(G)$ contains $\bigoplus_{\ensuremath{\mathbb{Z}}} G$ as a normal subgroup and fits into a short exact sequence
\begin{equation}
    1\rightarrow \bigoplus_{\ensuremath{\mathbb{Z}}} G\rightarrow \mathfrak{D}(G)\rightarrow G\rightarrow1,
\end{equation}
where the quotient map sends a sequence to its eventual constant value.

The following lemma is a consequence of the fact that countability, torsion-freeness and amenability pass to extensions and countable direct sums. 

\begin{lemma}
\label{lemma:E(G)-is-nice} 
    Whenever $G$ is countable, torsion-free and amenable, so is $\mathfrak{D}(G)$. 
\end{lemma} 

Consider the left shift automorphism of $\mathfrak{D}(G)$ given by 
\begin{equation}
    \sigma_G\colon \mathfrak{D}(G)\rightarrow \mathfrak{D}(G),
    \qquad
    (\sigma_G x)_k \coloneqq x_{k+1}.
\end{equation} 
We define the group $\mathfrak{C}(G)$ to be the semi-direct product 
\begin{equation}
    \mathfrak{C}(G)\coloneqq \mathfrak{D}(G)\rtimes_{\sigma_G}\ensuremath{\mathbb{Z}}.
\end{equation} 
We denote by \(t_G\) the generator of the acting group \(\ensuremath{\mathbb{Z}}\), with the convention that
\begin{equation}
    t_G x t_G^{-1}=\sigma_G(x), \quad x\in \mathfrak{D}(G).
\end{equation} 

In particular, we have $1\rightarrow \mathfrak{D}(G)\rightarrow \mathfrak{C}(G)\rightarrow \ensuremath{\mathbb{Z}}\rightarrow 1$. Combined with Lemma~\ref{lemma:E(G)-is-nice} and again that amenability, countability and
torsion-freeness pass to extensions, we get the following. 
\begin{lemma}
\label{lemma:S(G)-is-nice}
    Whenever $G$ is countable, torsion-free and amenable, so is $\mathfrak{C}(G)$. 
\end{lemma}

Note that there is a canonical coordinate-zero embedding
\begin{equation}
\label{eq:delta-g}
    \delta_G\colon G\hookrightarrow \mathfrak{D}(G)\leq \mathfrak{C}(G)
\end{equation}
defined by
\begin{equation}
    \delta_G(g)_k \coloneqq
    \begin{cases}
        g,& k=0,\\
        1,& k\neq 0.
    \end{cases}
\end{equation}

We now iterate the construction described above to obtain the desired acyclic group $\ensuremath{\mathfrak{B}}(G)$. 
Starting from a seed group \(G\), set \(G_0 \coloneqq G\), \(G_{n+1} \coloneqq \mathfrak{C}(G_n)\), with bonding map \(\delta_{G_n}\), and define
\begin{equation} 
\ensuremath{\mathfrak{B}}(G)\coloneqq \varinjlim_n (G_n, \delta_{G_{n}})=\bigcup_{n\ge0}G_n,
\end{equation}
where the second equality uses the injectivity of the bonding maps.
In particular, \(G\leq\ensuremath{\mathfrak{B}}(G)\).
Since countability, torsion-freeness and amenability are preserved by countable inductive limits, Lemma~\ref{lemma:S(G)-is-nice} gives the following proposition. 
\begin{prop}
\label{prop:B(G)-is-nice}
    The seed group $G$ embeds as a subgroup of $\ensuremath{\mathfrak{B}}(G)$.
    Moreover, $\ensuremath{\mathfrak{B}}(G)$ is countable, torsion-free and amenable whenever $G$ is.
\end{prop}

\subsection{Acyclicity of $\ensuremath{\mathfrak{B}}(G)$}

In this subsection, we show that $\ensuremath{\mathfrak{B}}(G)$ is always pseudo-mitotic and hence acyclic by Theorem~\ref{thm:pseudomitotic-acyclic}.
We begin by describing a group homomorphism that will facilitate pseudo-mitosis.

Let $G$ be a group.
Denote by $\rho_G$ the homomorphism given by 
\begin{equation}
    \rho_G\colon G\rightarrow \mathfrak{D}(G),\qquad\rho_{G}(g)_k \coloneqq \begin{cases}
        1,& k\leq 0,\\
        g,& k\geq 1.
    \end{cases}
\end{equation}
Since the image of $\rho_G$ has trivial $0$-th coordinate, it commutes with the coordinate-zero embedding $\delta_G$ of $G$ into $\mathfrak{D}(G)$, that is 
\begin{equation}
\label{eq:commutator}
    [\delta_G(g),\rho_G(g')]=1,\quad g,g'\in G. 
\end{equation}
Moreover, for all $g\in G$, we have the equality 
\begin{equation}
\label{eq:conjugation}
    \delta_G(g)\rho_G(g)=t_G\rho_G(g)t_G^{-1} 
\end{equation} 
inside $\mathfrak{C}(G)$, since conjugation by $t_G$ implements the left shift $\sigma_G$ on $\mathfrak{D}(G)$.

\begin{prop}
\label{prop:B-pseudomitotic}
For every seed group \(G\), the group \(\ensuremath{\mathfrak{B}}(G)\) is pseudo-mitotic.
\end{prop}

\begin{proof} 
Let \(H\leq \ensuremath{\mathfrak{B}}(G)\) be an arbitrary finitely generated subgroup. 
By finite generation, there exists some $n$ such that
\(H\leq G_n\). 
We will construct a pseudo-mitosis of \(H\) inside the next stage
\begin{equation}
    G_{n+1}=\mathfrak{C}(G_n)=\mathfrak{D}(G_n)\rtimes_{\sigma_{G_n}} \ensuremath{\mathbb{Z}}.
\end{equation} 
We identify $G_{n}$ with the coordinate-zero copy $\delta_{G_{n}}(G_{n})$ inside $G_{n+1}$. 
Next, define \mbox{$\psi_{0}, \psi_{1}\colon H \rightarrow G_{n+1}$} by 
\begin{equation} 
    \psi_{0}(h) = \delta_{G_{n}} (h) \rho_{G_{n}} (h) \quad \text{and} \quad \psi_{1}(h) = \rho_{G_{n}}(h). 
\end{equation} 
The map $\psi_{0}$ is a group homomorphism, as it is the product of two homomorphisms \(h\mapsto \delta_{G_n}(h)\) and \(h\mapsto \rho_{G_n}(h)\), whose images commute elementwise by \eqref{eq:commutator}. 
Take $g = t_{G_{n}} \in G_{n+1}$. It remains to check the three pseudo-mitosis conditions from Definition~\ref{def:pseudomitosis}.  
Condition \eqref{eq:pm1} holds by construction and condition \eqref{eq:pm2} holds by \eqref{eq:commutator}.
From \eqref{eq:conjugation}, we get 
\begin{equation}
    g\psi_1(h)g^{-1}
    =t_{G_n}\rho_{G_n}(h)t_{G_n}^{-1}
    =h\rho_{G_n}(h) 
    =\psi_0(h), 
\end{equation}
which verifies condition \eqref{eq:pm3}. Since \(H\) was  arbitrary, \(\ensuremath{\mathfrak{B}}(G)\) is pseudo-mitotic. 
\end{proof}

Acyclicity of $\ensuremath{\mathfrak{B}}(G)$ follows from the Varadarajan-Berrick theorem (Theorem~\ref{thm:pseudomitotic-acyclic}).

\begin{cor}
\label{cor:B-acyclic}
    For every seed group \(G\), the group \(\ensuremath{\mathfrak{B}}(G)\) is acyclic.
\end{cor} 

\begin{proof}[Proof of Theorem~\ref{introthm:embedings-into-acyclic}] 
    Let $G$ be an amenable, countable, torsion-free group. 
    By Proposition~\ref{prop:B(G)-is-nice}, the group $\ensuremath{\mathfrak{B}}(G)$ is amenable, countable torsion-free, and contains $G$ as a subgroup.
    By Corollary~\ref{cor:B-acyclic}, $\ensuremath{\mathfrak{B}}(G)$ is acyclic.
\end{proof}

\section{An uncountable family of examples} 
\label{section:uncountable_family} 

In this section we construct a continuum family of countable, torsion-free, amenable and 
acyclic groups and show that, for every fixed countable group, a continuum
subfamily remains distinguishable after taking free products. The proofs use the following isomorphism invariant:
whether the group contains a non-trivial element that is divisible by arbitrarily high powers of
a given prime.
We denote by $\ensuremath{\mathcal{P}}$ the set of all primes. 

\begin{defn}
\label{def:D}
Let \(G\) be a group and let \(p\) be a prime.  
A non-identity element \(g\in G\) is called \emph{infinitely \(p\)-divisible} if, for every \(r\geq 1\), there exists \(x_r\in G\) such that
\begin{equation}
    x_r^{p^r}=g.
\end{equation}
Define
\begin{equation}
    D(G) \coloneqq \{p\in \ensuremath{\mathcal{P}} \ | \  G \text{ contains an infinitely } p\text{-divisible element}\}.
\end{equation} 
\end{defn}
The set \(D(G)\) is an isomorphism invariant of \(G\). 
We are of course interested in understanding $D(\ensuremath{\mathfrak{B}}(G))$ in terms of $D(G)$, as shown in the following proposition. 

\begin{prop}
\label{prop:D-BA}
Let $G$ be any group. 
Then 
\(
    D(\ensuremath{\mathfrak{B}}(G))=D(G).
\)
\end{prop}

\begin{proof}
Set \(G_0=G\) and \(G_{n+1}=\mathfrak{C}(G_n)\), so that
\(\ensuremath{\mathfrak{B}}(G)=\varinjlim G_n\). 
Since the bonding maps are injective, we identify each $G_{n}$ with its image in $\ensuremath{\mathfrak{B}}(G)$. 
We first observe that \(D(\mathfrak{C}(H))=D(H)\) for every group \(H\).  
The inclusion $D(H) \subseteq D(\mathfrak{C}(H))$ follows from the coordinate-zero embedding $\delta_{H}$.  
Conversely, suppose that \(u\in \mathfrak{C}(H)=\mathfrak{D}(H)\rtimes \ensuremath{\mathbb{Z}}\) is a non-trivial infinitely \(p\)-divisible element. 
Since non-trivial elements in $\ensuremath{\mathbb{Z}}$ are not infinitely divisible, $u$ and all its \(p^r\)-roots belong to $\mathfrak{D}(H)$ for any $r\geq 1$.  
Evaluating at a coordinate on which \(u\) is non-trivial gives an infinitely \(p\)-divisible element of \(H\). 

Next we show that roots are reflected across the bonding maps. For any group $H$, an element $h\in H$ and $m\geq 1$, suppose that $\delta_{H}(h) = x^{m}$ in $\mathfrak{C}(H)$. 
By the definition of the bonding map $\delta_{H}$, it follows that $x\in \mathfrak{D}(H)$, and evaluating at coordinate-zero, we have $x_{0}^{m} = h$. 
Thus, whenever \(\delta_H(h)\) has an $m$-th root in $\mathfrak{C}(H)$, the element $h$ has an $m$-th root in $H$.

The inclusion $D(G) \subseteq D(\ensuremath{\mathfrak{B}}(G))$ follows from the embedding $G\hookrightarrow \ensuremath{\mathfrak{B}}(G)$. 
For the reverse inclusion, let $ p\in D(\ensuremath{\mathfrak{B}}(G))$, and choose an infinitely $p$-divisible element \(g\in\ensuremath{\mathfrak{B}}(G)\). 
In particular, $g\in G_{n}$ for some $n\in \mathbb{N}$. 
For each $r\geq 1$, there exists $x_{r}\in \ensuremath{\mathfrak{B}}(G)$ such that $x_{r}^{p^{r}} = g$, and choose some $n_{r}\geq n$ such that $x_{r}\in G_{n_{r}}$. 
Applying the argument in the paragraph above shows that a $p^{r}$-root of $g$ can be found in $G_{n}$. 
Consequently, $g$ is infinitely $p$-divisible in $G_{n}$, and hence $p\in D(G_{n}) = D(G)$ by repeatedly applying $D(\mathfrak{C}(H))=D(H)$. Thus, $D(\mathfrak B(G))=D(G)$. 
\end{proof}

We now use a construction due to Hall (\cite{hall}) to produce uncountably many finitely generated amenable groups with distinct $D$ invariants. 
For a group $G$, we denote by $Z(G)$ the \emph{centre} of the group, that is, 
\begin{equation} 
Z(G) = \{g\in G \ | \ gk = kg \text{ for every }k\in G\}. 
\end{equation} 

\begin{lemma}
\label{lem:solvable-seeds}
For every subset \(S\subseteq\ensuremath{\mathcal{P}}\), there is a two-generated,
solvable, torsion-free group \(M_S\) such that \(D(M_S)=S\).
\end{lemma}

\begin{proof} 
Take $C_S= (\ensuremath{\mathbb{Z}}[S^{-1}], +)$ obtained from $\ensuremath{\mathbb{Z}}$ by allowing division by the primes in $S$. 
Note that $C_S$ is a torsion-free group and $D(C_{S}) = S$ by definition. 
By \cite[Lemma~4.1]{caprace2024simple}, following a classical construction of Hall, there is a two-generated, three-step
solvable, torsion-free group \(M_S\) such that
\(Z(M_S)\cong C_S\), \(M_S/Z(M_S)\cong\ensuremath{\mathbb{Z}}\wr\ensuremath{\mathbb{Z}}\).

Let \(W=\ensuremath{\mathbb{Z}}\wr\ensuremath{\mathbb{Z}}=(\bigoplus_{\ensuremath{\mathbb{Z}}}\ensuremath{\mathbb{Z}})\rtimes \ensuremath{\mathbb{Z}}\).  Then
\(D(W)=\varnothing\).  Indeed, projection onto the acting copy of \(\ensuremath{\mathbb{Z}}\)
places an infinitely \(p\)-divisible element, together with all of its roots,
in \(\bigoplus_{\ensuremath{\mathbb{Z}}}\ensuremath{\mathbb{Z}}\), where coordinatewise divisibility forces it to be
trivial.

Since \(D(C_S)=S\), it remains only to show that no new divisibility occurs
in \(M_S\).  If \(g\in M_S\) is infinitely \(p\)-divisible, then its image in
\(W\) is trivial, so \(g\in Z(M_S)\).  Every root of \(g\) also belongs to
\(Z(M_S)\), since \(W\) is torsion-free.  Thus
\(p\in D(Z(M_S))=S\).  The reverse inclusion follows from
\(Z(M_S)\cong C_S\).
\end{proof}

We get the following immediate result. 
Note that for every group $G$, $\ensuremath{\mathfrak{B}}(G)$ is non-cyclic; in fact, it is infinitely generated, as a consequence of pseudo-mitosis (\cite[Theorem 1.8]{Berrick1989}).

\begin{cor}
\label{cor:acyclic-separating-family} 
For every subset $S\subseteq \mathcal{P}$, the group \(\ensuremath{\mathfrak{B}}(M_S)\) is countable, amenable, acyclic, torsion-free, non-cyclic and satisfies \(D(\ensuremath{\mathfrak{B}}(M_S))=S\).  
In particular, groups in $\{\ensuremath{\mathfrak{B}}(M_S) \ | \ S\subseteq \mathcal{P}\}$ are pairwise non-isomorphic. 
\end{cor}

We are now ready to prove the main mechanism for producing uncountably many groups which distinguish the free products. 

\begin{thm}
\label{prop:free-product-separation}
Let \(\Gamma\) be a countable group.  There is a set
\(\mathcal S_\Gamma\subseteq 2^{\ensuremath{\mathcal{P}}}\) of cardinality \(2^{\aleph_0}\) such
that \(M_S\) does not embed into \(\Gamma\) for every
\(S\in\mathcal S_\Gamma\), and the groups \(\Gamma*\ensuremath{\mathfrak{B}}(M_S)\),
\(S\in\mathcal S_\Gamma\), are pairwise non-isomorphic and are not
isomorphic to \(\Gamma\).
\end{thm} 

\begin{proof} 
A countable group has only countably many two-generated subgroups.  Since the
\(M_S\) are two-generated and pairwise non-isomorphic, only countably many of
them can embed into \(\Gamma\).  Therefore
\begin{equation}
    \mathcal S_\Gamma
    =\{S\subseteq\ensuremath{\mathcal{P}} \ | \ M_S\text{ does not embed into }\Gamma\}
\end{equation}
has cardinality \(2^{\aleph_0}\). 

Fix distinct \(S,T\in\mathcal S_\Gamma\), and suppose that
\(\theta\colon\Gamma*\ensuremath{\mathfrak{B}}(M_S)\to\Gamma*\ensuremath{\mathfrak{B}}(M_T)\) is an isomorphism.  The subgroup
\(\theta(\ensuremath{\mathfrak{B}}(M_S))\) is non-cyclic, torsion-free and amenable, so
Lemma~\ref{lem:amenable-subgroups-free-product}(i) places it inside a conjugate
of \(\Gamma\) or of \(\ensuremath{\mathfrak{B}}(M_T)\).  The first possibility would give an embedding
\(M_S\hookrightarrow\Gamma\), contrary to the definition of
\(\mathcal S_\Gamma\).  Hence
\(\theta(\ensuremath{\mathfrak{B}}(M_S))\leq g\ensuremath{\mathfrak{B}}(M_T)g^{-1}\) for some \(g\in\Gamma*\ensuremath{\mathfrak{B}}(M_T)\).

By Lemma~\ref{lem:amenable-subgroups-free-product}(ii), \(\ensuremath{\mathfrak{B}}(M_S)\) is maximal
amenable in \(\Gamma*\ensuremath{\mathfrak{B}}(M_S)\); hence \(\theta(\ensuremath{\mathfrak{B}}(M_S))\) is maximal amenable in the
target.  Since \(g\ensuremath{\mathfrak{B}}(M_T)g^{-1}\) is amenable, the preceding inclusion is an
equality.  Thus \(\ensuremath{\mathfrak{B}}(M_S)\cong \ensuremath{\mathfrak{B}}(M_T)\), and consequently
\begin{equation}
S=D(\ensuremath{\mathfrak{B}}(M_S))=D(\ensuremath{\mathfrak{B}}(M_T))=T,
\end{equation}
a contradiction.

Finally, an isomorphism \(\Gamma*\ensuremath{\mathfrak{B}}(M_S)\cong\Gamma\) would restrict to an
embedding \(M_S\hookrightarrow\Gamma\), which is impossible.
\end{proof}

\section{Necessary classification results} 
\label{section:classification} 
As stated in the introduction, our proof relies heavily on the recent breakthrough isomorphism theorem of Hirshberg and Phillips (\cite{hp}) for infinite reduced free products of \ensuremath{C^{*}}-algebras. 
To the best of our understanding, their result (\cite[Corollary~3.4]{hp}) does not cover the \ensuremath{C^{*}}-algebras of the torsion-free, amenable, acyclic groups constructed in Section~\ref{sec:construction-of-groups}.
The required extension will follow from the uniqueness theorem below, which is a special case of the main result in the upcoming paper \cite{classification2}. 
We kindly thank the authors of \cite{classification2} for pointing out the result and letting us include it here. 

\subsection{Classification of morphisms into monotracial \ensuremath{C^{*}}-algebras} 
We include a self-contained proof of Theorem~\ref{thm:classification} for completeness, by using the existing trace-kernel machinery and the recent breakthrough $KK$-uniqueness result by Szab\'o. We refer to the classification paper \cite{classification} for basics on Cuntz pairs, $KK$-theory and the UCT, (unital) absorption of full maps, etc.

Since we are in the unique trace case, it is more convenient to work with ultrapowers rather than sequence algebras. 
Fix a free ultrafilter $\omega$ on $\mathbb{N}$, and let $B$ be a unital, simple, non-elementary $C^*$-algebra which has strict comparison with respect to a unique trace $\tau_{B}$. 
By \cite[Theorem 1.2]{ngrobert}, $B_{\omega}$ has a unique trace $\tau_{B_{\omega}}$, given by 
\begin{equation} 
\tau_{B_{\omega}}((b_n)_{n}) \coloneqq \lim_{n\to \omega} \tau_{B}(b_n).
\end{equation}
Then the \emph{trace--kernel ideal} of $B_\omega$ is given by 
\begin{equation}
J_B\coloneqq \{x\in B_\omega \ | \ \tau_{B_{\omega}}(x^*x)=0\}, 
\end{equation} 
and it gives rise to the \emph{trace--kernel extension}
\begin{equation} \label{Trace_kernel}
0\longrightarrow J_B\stackrel{j_B}{\longrightarrow} B_\omega\stackrel{q_B}{\longrightarrow} B^\omega \longrightarrow 0, 
\end{equation} 
where $B^\omega\coloneqq B_\omega/J_B$ is the \emph{tracial ultrapower} of $B$. The trace $\tau_{B_\omega}$ descends to a faithful trace $\tau_{B^\omega}$ on $B^\omega$, that is, $\tau_{B_{\omega}} = \tau_{B^{\omega}}\circ q_{B}$. 
Since $\tau_{B}$ is the unique trace on $B$, the von Neumann algebra $\pi_{\tau_B}(B)''$ is a II$_{1}$ factor (\cite[Theorem 6.7.4]{Dixmier}). Moreover, $B^{\omega}$ canonically identifies with the tracial von Neumann algebra ultraproduct 
$(\pi_{\tau_{B}}(B)'')_{\omega}$, which is again a II$_{1}$ factor (\cite[Remark 4.7]{KR_SI}).

\begin{thm}
\label{thm:classification} 
Let $A$ be a unital, separable, nuclear \ensuremath{C^{*}}-algebra satisfying the UCT, with $K_0(A)\cong \ensuremath{\mathbb{Z}}\cdot [1]_0 $ and $K_{1}(A) = 0$. 
Let $B$ be a unital, simple, exact, non-elementary \ensuremath{C^{*}}-algebra which has strict comparison with respect to a unique trace $\tau_{B}$, and $K_{*}(B) = (\ensuremath{\mathbb{Z}} \cdot [1_{B}]_{0},0)$. 

If $\phi,\psi\colon A\rightarrow B$ are unital embeddings satisfying $\tau_{B} \circ \phi = \tau_{B} \circ \psi$, 
then $\phi$ and $\psi$ are approximately unitarily equivalent. 
\end{thm} 

\begin{proof}[Proof of Theorem~\ref{thm:classification}]
Let $\iota_{B}\colon B \rightarrow B_{\omega}$ denote the diagonal inclusion, and set 
$\Phi = \iota_{B} \circ \phi$, $\Psi = \iota_{B} \circ \psi$, which are unital embeddings of $A$ into $B_{\omega}$. 

We first classify the induced maps into $B^{\omega}$. 
Since $B$ is simple, it follows that $\Phi$ and $\Psi$ are unital and full, which implies that $q_{B} \circ \Phi$ and $q_{B} \circ \Psi$ are also unital and full. 
In particular, $q_{B} \circ \Phi$ and $q_{B} \circ \Psi$ are unital embeddings into the II$_{1}$ factor $B^{\omega}$. 
The equality $\tau_{B} \circ \phi = \tau_{B} \circ \psi$ implies that 
\begin{equation} 
    \tau_{B^{\omega}} \circ q_{B} \circ \Phi = \tau_{B_{\omega}} \circ \Phi = \tau_{B} \circ \phi = \tau_{B} \circ \psi = \tau_{B_{\omega}} \circ \Psi = \tau_{B^{\omega}} \circ q_{B} \circ \Psi, 
\end{equation} 
where the first and last equalities follow from $\tau_{B_{\omega}} = \tau_{B^{\omega}} \circ q_{B}$. 

By uniqueness of nuclear maps into finite factors, a consequence of Connes' theorem (\cite[Proposition 1.1]{Schafhauser.Annals}), $q_B\circ\Phi,q_B\circ\Psi\colon A\to B^\omega$ are approximately unitarily equivalent in the $2$-norm. 
Separability of $A$ and a reindexing argument in $B^\omega$ then give unitary equivalence. 
Every unitary in a II$_1$ factor is homotopic to the identity and therefore lifts to a unitary in $B_\omega$. 
After replacing $\Psi$ by a unitary conjugate, we may assume that $q_B\circ\Phi=q_B\circ\Psi$.  

The following properties of $J_{B}$, $B_{\omega}$ and $B^{\omega}$ are consequences of our assumptions: 
\begin{enumerate}[label=(\roman*)] 
\item 
Since $B$ is simple and has strict comparison with respect to $\tau_{B}$, the norm ultrapower $B_{\omega}$ has strict comparison with respect to $\tau_{B_{\omega}}$ (\cite[Lemma 1.23]{BBSTWW}). 
Moreover, with the additional assumption $K_{1}(B) = 0$, we have $K_{1}(B_{\omega}) = 0$ by \cite[Corollary 4.2]{SriSchafhauser}; 
\item 
The tracial ultrapower $B^{\omega}$ is a II$_{1}$ factor, thus $K_{0}(B^{\omega}) \cong \mathbb{R}$ is divisible and $K_{1}(B^{\omega}) = 0$; 
\item 
The argument of \cite[Lemma 6.11]{classification} (using strict comparison of $B_{\omega}$ instead of $\ensuremath{\mathcal{Z}}$-stability of $B$) implies that $J_{B}$ is separably stable, which means that every separable $C^*$-subalgebra is contained in a stable and separable $C^*$-subalgebra. 
Lastly, since $\mathrm{Cu}(B)$ is almost unperforated, so is $\mathrm{Cu}(J_{B})$ as the property passes to ultrapowers and ideals. 
\end{enumerate} 
Moreover, since $K_{0}(A) = \mathbb{Z}\cdot [1]_{0}$ and $\phi, \psi$ are unital, it follows 
that $K_{0}(\phi) = K_{0}(\psi)$, and hence $K_{0}(\Phi) = K_{0}(\Psi)$. 
Simplicity of $B$ implies that $\Phi$ and $\Psi$ are full maps. 

As is standard, we proceed with a separabilization argument before computing $KK$-theory. 
Properties including $K_{0}$ being divisible, trivial $K_{1}$, almost unperforation of the  Cuntz semigroup are separably inheritable (\cite[Proposition 4.1]{new-TWW}). We now separabilize using the same approach as \cite[Section 1.2]{Schafhauser.Annals} to obtain a separable unital subextension 
\begin{equation}
\begin{tikzcd} 
0\arrow[r]&J\arrow[r,"j"]\arrow[d]&E\arrow[r,"q"]\arrow[d]&D\arrow[r]\arrow[d]&0\\
0\arrow[r]&J_B\arrow[r,"j_B"]&B_{\omega}\arrow[r,"q_B"]&B^\omega\arrow[r]&0
\end{tikzcd}
\end{equation} 
such that the following hold:
\begin{enumerate}[label=(\roman*)]
\item \label{condition_sep_full} $\Phi(A)\cup \Psi(A)\subseteq E$, and the corestrictions $\Phi|^E,\Psi|^E\colon A\to E$ of $\Phi$ and $\Psi$ to $E$ are unital, full, and satisfy $K_{0}(\Phi|^E)=K_{0}(\Psi|^E)$; 
\item $K_{1}(D)$, $K_{1}(E)$ vanish, and $K_{0}(D)$ is divisible; 
\item \label{condition_sep_J}$\mathrm{Cu}(J)$ is almost unperforated, and $J$ is stable. 
\end{enumerate} 

Let $\lambda\colon E\to \mathcal M(J)$ be the canonical unital $^*$-homomorphism. 
Since $q_B\circ \Phi = q_{B}\circ \Psi$, we have $q\circ \Phi|^E = q\circ \Psi|^E$ and hence 
\begin{equation}
\left(\lambda\circ\Phi|^E,\lambda\circ\Psi|^E \right)\!\colon A\rightrightarrows \ensuremath{\mathcal{M}}(J)\rhd J
\end{equation}
is a Cuntz pair inducing a class $\kappa\in KK(A,J)$, which satisfies 
\begin{equation} 
KK(A,j)(\kappa)=[\Phi|^E]_{KK(A,E)}-[\Psi|^E]_{KK(A,E)}
\end{equation}
by \cite[Proposition 2.1]{Schafhauser.Annals}. 
Since $A$ satisfies the UCT, triviality of $K_{1}(A)$ and $K_{1}(E)$ implies that there is an isomorphism $KK(A, E) \cong \mathrm{Hom}(K_{0}(A), K_{0}(E))$, and hence $[\Phi|^{E}]_{KK(A, E)} = [\Psi|^{E}]_{KK(A, E)}$ as $K_{0}(\Phi|^E)=K_{0}(\Psi|^E)$. 
Therefore $KK(A,j)(\kappa) = 0$. 
To prove $\kappa = 0$, it suffices to prove $KK(A, j)$ is injective. 
The relevant part of the six-term sequence in \(KK\)-theory  is
\begin{equation} 
KK^1(A,D) = KK(A, SD)\xrightarrow{\partial}KK(A,J) 
   \xrightarrow{KK(A, j)}KK(A,E). 
\end{equation} 
where $SD$ denotes the suspension of $D$. 
It therefore suffices to prove $KK(A, SD) = 0$. 
Indeed, again by the UCT of $A$, $K_{1}(A) = 0$ and $K_{1}(D) = 0$, there is an isomorphism
\begin{equation} 
\mathrm{Ext}^1(K_0(A), K_{1}(SD)) \cong KK(A,SD) 
\end{equation} 
which is trivial as $K_{1}(SD) \cong K_{0}(D)$ is divisible and is therefore injective as an abelian group. 
Thus, 
\begin{equation} 
\kappa = \left[ \lambda\circ\Phi|^E,\lambda\circ\Psi|^E \right]_{KK(A, J)} = 0. 
\end{equation} 

Almost unperforation of $\mathrm{Cu}(J)$ implies that $J$ has the corona factorisation property by the argument of \cite[Proposition 4.12]{classification} (which only uses almost unperforation as a consequence of $\ensuremath{\mathcal{Z}}$-stability). 
Since $J$ is stable and has the corona factorisation property, and $\lambda \circ \Phi|^{E}, \lambda \circ \Psi|^{E}$ are full, it follows that  $\lambda\circ \Phi|^E$ and $\lambda\circ \Psi|^E$ are unitally absorbing by Elliott-Kucerovsky's generalized Weyl-von Neumann theorem (\cite{EK.PJM}). 
Combining Szab\'o's $K_{1}$-injectivity result (\cite[Corollary 2.9]{gabor}), with the standard $KK$-uniqueness theorem for unitally absorbing maps (see \cite[Theorem 1.12]{huawhite} for instance), there exists a sequence of unitaries $(u_n)_{n}$ in the minimal unitization $J^{\dagger}$ of $J$ such that 
\begin{equation}
\|u_n(\lambda\circ  \Phi|^E(a))u_n^*-\lambda\circ  \Psi|^E(a)\|\to 0,\quad a\in A. 
\end{equation}
Since $\lambda$ restricts to the identity map on $J^\dagger$, it follows that  
\begin{equation} 
\|u_n  \Phi|^E(a)u_n^*- \Psi|^E(a)\|\to 0,\quad a\in A. 
\end{equation} 
Thus the unital $^*$-homomorphisms $\Phi, \Psi\colon A\rightarrow B_{\omega}$ are approximately unitarily equivalent, which implies that $\phi$ and $\psi$ are approximately unitarily equivalent. 
\end{proof}

\subsection{The generalized Hirshberg-Phillips theorem}
This subsection contains the proof of Theorem~\ref{introthm:hp-upgrade}, the generalization of the isomorphism theorem of Hirshberg and Phillips (\cite[Theorem~3.1]{hp}). 
A direct corollary of Theorem~\ref{thm:classification} is the following uniqueness theorem of maps into reduced free products. 

\begin{cor}
\label{cor:replacement-for-3.4}
    Let $A$ be a unital, separable, nuclear \ensuremath{C^{*}}-algebra satisfying the UCT, with $K_{*}(A) = (\mathbb{Z}\cdot [1]_{0}, 0)$, and suppose $A$ has a faithful trace $\tau_A$. 
    Let $\theta\colon A\rightarrow\ensuremath{\mathcal{Z}}$ be an embedding satisfying $\tau_A=\tau_\ensuremath{\mathcal{Z}}\circ \theta$, where $\tau_\ensuremath{\mathcal{Z}}$ is the unique trace on $\ensuremath{\mathcal{Z}}$. 
    Let $A*_r\ensuremath{\mathcal{Z}}$ denote the reduced free product with respect to the traces $\tau_A$ and $\tau_\ensuremath{\mathcal{Z}}$, and let 
    \begin{equation}
        \begin{aligned}
            \kappa_1&\colon A\rightarrow A*_r\ensuremath{\mathcal{Z}},\\
            \kappa_2&\colon \ensuremath{\mathcal{Z}}\rightarrow A*_r\ensuremath{\mathcal{Z}}
        \end{aligned}
    \end{equation}
    be the canonical embeddings.
    Then $\kappa_1$ is approximately unitarily equivalent to $\kappa_2\circ\theta$.
\end{cor}

\begin{proof}
    By \cite[Corollary~5.3 and Theorem~3.1]{robert}, $A*_r\ensuremath{\mathcal{Z}}$ is selfless and in particular is simple and has strict comparison with respect to its unique trace $\tau$. 
    Since $\tau \circ \kappa_{1} = \tau \circ \kappa_2\circ\theta = \tau_{A}$ and $A*_r\ensuremath{\mathcal{Z}} $ has a unique trace, both maps agree on the tracial data. 
    Moreover, $A *_{r} \ensuremath{\mathcal{Z}}$ is exact by \cite[Corollary~4.3]{dykemashlyakhtenko}. 
    A combination of \cite[Theorem~6.4]{thomsen03} and \cite[Theorem~1.1]{hasegawa16} shows that $K_{*}(A*_{r}\ensuremath{\mathcal{Z}}) = (\mathbb{Z}\cdot [1]_{0}, 0)$. 
    Theorem~\ref{thm:classification} therefore applies and shows that the two maps are approximately unitarily equivalent.
\end{proof}

With Corollary~\ref{cor:replacement-for-3.4} established, the proof of  \cite[Theorem~3.1]{hp} can be modified to remove the requirement that $A$ is a limit of 1-dimensional NCCW complexes.

\begin{thm}[{A combination of \cite{hp}, \cite{classification} and \cite{gabor}}]
    \label{thm:hirshberg-phillips}
    Let $\{C_i\}_{i\in \ensuremath{\mathbb{N}}}$, \mbox{$C_i\neq \ensuremath{\mathbb{C}}$} be  unital, separable \ensuremath{C^{*}}-algebras equipped with faithful traces $\tau_{C_i}\in T(C_i)$ and such that each $C_i$ contains a self-adjoint element with diffuse spectral measure (induced by $\tau_{C_i}$).
    Let $A$ be a unital, separable, nuclear \ensuremath{C^{*}}-algebra satisfying the UCT and equipped with a faithful trace $\tau_A$. 
    Suppose  that $K_*(A)= (\ensuremath{\mathbb{Z}}\cdot[1]_0,0)$.
    Then
    \begin{equation}
        A*_r \mathop{\scalebox{1.7}{$\;*$}\,\,}\limits_{i\in \ensuremath{\mathbb{N}}} C_i \cong \mathop{\scalebox{1.7}{$\;*$}\,\,}\limits _{i\in\ensuremath{\mathbb{N}}} C_i,
    \end{equation}
    where the free products are taken with respect to the traces $\tau_A$ and $\tau_{C_i}$.
\end{thm}

The proof is a small modification of \cite[Theorem~3.1]{hp} so we only point out the necessary adjustments.

\begin{proof}
    First note that $C_i$ containing a self-adjoint element with diffuse spectral measure implies that $\operatorname{At}(C_i,\tau_{C_i})=0$ for all $i\in \ensuremath{\mathbb{N}}$.
    The quantity $\operatorname{At}(C_i,\tau_{C_i})$ is defined in \cite[Definition~2.5]{hp} and measures atoms in $C_i$ with respect to the measure induced by $\tau_{C_i}$.
    Thus the applications of \cite[Lemma~2.6]{hp} remain valid.
    
    The assumption in \cite[Proposition~2.13]{hp} that $A$ is an inductive limit of \mbox{1-dimensional} NCCW complexes is redundant, see \cite[Corollary~C]{classification}.
    Hence there is an embedding $\theta\colon A\rightarrow \ensuremath{\mathcal{Z}}$ such that $\tau_{A} = \tau_{\ensuremath{\mathcal{Z}}} \circ \theta$, 
    where $\tau_A$ is the distinguished trace on $A$. Corollary~\ref{cor:replacement-for-3.4} therefore applies and shows that $\kappa_1$ and $\kappa_2\circ\theta$ are approximately unitarily equivalent.
    
    The proof of \cite[Theorem~3.1]{hp} can then proceed word-for-word with $A$ (with the trace $\tau_A$) and $C_i$ (with the trace $\tau_{C_i}$) after substituting  \cite[Corollary~3.4]{hp} with Corollary~\ref{cor:replacement-for-3.4} (and keeping in mind that the conclusion of \cite[Proposition~2.13]{hp} remains valid for $A$).
\end{proof}

\section{Counterexamples to \ensuremath{C^{*}}-superrigidity} 
\label{sec:counterexample} 

This section derives the main \ensuremath{C^{*}}-superrigidity result, Theorem~\ref{introthm:general-nonrigidity}, which provides the first examples of torsion-free groups which are not \ensuremath{C^{*}}-superrigid.
The examples arise as free products of $G$ with an infinite free product, where $G$ is discrete, countable, torsion-free, amenable and acyclic. 

The following is an adaptation of the Hirshberg-Phillips theorem (Theorem~\ref{thm:hirshberg-phillips}) to group \ensuremath{C^{*}}-algebras.
\begin{thm}
\label{thm:iso}
    Let $G$ be an amenable, countable, torsion-free and acyclic group, and let $\{H_i\}_{i\in \ensuremath{\mathbb{N}}}$ be non-trivial countable groups.
    Then 
    \begin{equation}
        \ensuremath{C^{*}_r}\left(G*\mathop{\scalebox{1.7}{$\;*$}\,\,}\limits _{i\in\ensuremath{\mathbb{N}}}H_i\right)\cong \ensuremath{C^{*}_r}\left(\mathop{\scalebox{1.7}{$\;*$}\,\,}\limits _{i\in\ensuremath{\mathbb{N}}}H_i\right).
    \end{equation}
\end{thm}

\begin{proof}
    Note that $\ensuremath{C^{*}_r}(G)$ is a unital, separable, nuclear \ensuremath{C^{*}}-algebra, and is equipped with its canonical faithful trace $\tau_G$. 
    By \cite[Proposition~10.7]{tu}, it satisfies the UCT. 
    Moreover, by Lemma~\ref{lem:baum-connes}, we have
    \begin{align}
            K_0(\ensuremath{C^{*}_r}(G))= \ensuremath{\mathbb{Z}}\cdot [1]_0 \quad \text{and} \quad 
            K_1(\ensuremath{C^{*}_r}(G))= 0.
    \end{align} 

    For non-trivial groups $H_{1}$ and $H_{2}$, we choose non-identity elements $a\in H_{1}$ and $b\in H_{2}$. The reduced word $(ab)^{n}$ is nontrivial in $H_1* H_2$ for every $n\neq 0$, so $ab$ has infinite order. 
    Take $u_{ab}$ to be the canonical unitary in $C_{r}^{*}(H_1 * H_2)$, then $\displaystyle{\tau_{H_1 * H_2}\left(u_{ab}^{n}\right) = 0}$ for every $n\neq 0$. 
    Hence the induced spectral measure on $u_{ab}$ is the Haar measure on $\mathbb{T}$. Therefore, $(u_{ab} + u_{ab}^{*}) / 2$ is a self-adjoint element in $C_{r}^{*}(H_1 * H_2)$ with diffuse spectral measure (induced by $\tau_{H_1 * H_2}$). 
    By replacing $H_i$ with $H_{2i-1}*H_{2i}$, we may therefore assume that $\ensuremath{C^{*}_r}(H_i)$ contains a self-adjoint element with diffuse spectral measure (induced by $\tau_{H_i}$) for every $i\in \mathbb{N}$. 

    An application of Theorem~\ref{thm:hirshberg-phillips} to $A=\ensuremath{C^{*}_r}(G)$ with its canonical trace $\tau_G$ and $C_i=\ensuremath{C^{*}_r}(H_i)$ with their canonical traces $\tau_{H_i}$, shows that there is an isomorphism
    \begin{equation}
        \ensuremath{C^{*}_r}(G)*_r \mathop{\scalebox{1.7}{$\;*$}\,\,}\limits _{i\in\ensuremath{\mathbb{N}}}\ensuremath{C^{*}_r}(H_i)\cong \mathop{\scalebox{1.7}{$\;*$}\,\,}\limits _{i\in\ensuremath{\mathbb{N}}}\ensuremath{C^{*}_r}(H_i),
    \end{equation} 
    where the free products are taken with respect to the traces $\tau_G$ and $\tau_{H_i}$.
    The result then follows from Lemma~\ref{lemma:free-product-correspondnce}.
\end{proof}

Using the construction in Section~\ref{section:uncountable_family}, we get the following, which proves Theorem~\ref{introthm:general-nonrigidity}.

\begin{thm}
\label{thm:general-nonrigidity}
Let \(G\cong\ast_{i\in\ensuremath{\mathbb{N}}}G_i\), where every \(G_i\) is non-trivial and countable.  Then
there are continuum many pairwise non-isomorphic groups \(K\), none isomorphic
to \(G\), such that
\begin{equation}
    \ensuremath{C^{*}_r}(K)\cong\ensuremath{C^{*}_r}(G).
\end{equation}
If \(G\) is torsion-free, the groups \(K\) may be chosen torsion-free.
\end{thm} 

\begin{proof} 
Apply Theorem~\ref{prop:free-product-separation} with \(\Gamma=G\).  For every
\(S\in\mathcal S_G\), put \(K_S=G*\ensuremath{\mathfrak{B}}(M_S)\).  Each \(\ensuremath{\mathfrak{B}}(M_S)\) is torsion-free,
amenable and acyclic.  Applying Theorem~\ref{thm:iso} with the amenable factor
\(\ensuremath{\mathfrak{B}}(M_S)\) and with \(H_i=G_i\) gives 
\begin{equation}
    \ensuremath{C^{*}_r}(K_S)\cong\ensuremath{C^{*}_r}(G), 
    \quad S\in\mathcal S_G.
\end{equation}
Theorem~\ref{prop:free-product-separation} shows that the groups \(K_S\)
are pairwise non-isomorphic and are not isomorphic to \(G\).  If \(G\) is
torsion-free, then so is every \(K_S\), by
Lemma~\ref{lemma:torsion-free-products}.
\end{proof}

Taking \(G_i=\ensuremath{\mathbb{Z}}\) for every \(i\in \mathbb{N}\) in the preceding theorem gives a negative answer to the open question whether $F_{\infty}$ is \ensuremath{C^{*}}-superrigid. 

\begin{cor}
\label{cor:F_infty}
    The free group $F_{\infty}$ on infinitely many generators is not \ensuremath{C^{*}}-superrigid.
\end{cor} 

\section{Realizing the Jiang--Su algebra as a crossed product} 
\label{sec:model-for-Z} 

In this section we construct an action $G\curvearrowright X$ whose crossed product is isomorphic to the Jiang--Su algebra $\ensuremath{\mathcal{Z}}$ (Theorem~\ref{introthm:Z_realize}). 
We first recall the following results of Deeley, Putnam, and Strung (\cite{DPS-jiang-su}) that provide the base space and the starting point of the action. A compact, Hausdorff space $X$ is called \emph{point-like} if $C(X)$ is $KK$-equivalent to $\mathbb{C}$. 

\begin{prop}[{\cite[Corollary~1.12, Theorems~1.16--1.17]{DPS-jiang-su}}] 
\label{space_Z} 
There is an infinite, connected, compact, point-like metric space $Z$ with finite covering dimension. Moreover, there is a free, minimal, uniquely ergodic homeomorphism $\zeta$ of $Z$. 
\end{prop} 

\begin{proof} 
Fix an odd integer $d>1$ and a minimal, uniquely ergodic diffeomorphism $\varphi$ on $S^{d}$. 
The construction in \cite{DPS-jiang-su} yields a topological space $Z$ and a homeomorphism $\zeta$ of $Z$. 
The first statement follows from \cite[Corollary~1.12]{DPS-jiang-su}, and minimality is shown in \cite[Theorem 1.16]{DPS-jiang-su}. 
Since $\zeta$ is a minimal homeomorphism on an infinite space, every orbit is infinite, hence $\zeta$ is free. 
Lastly, $\zeta$ is uniquely ergodic by \cite[Theorem 1.17]{DPS-jiang-su} and unique ergodicity of $\varphi$. 
\end{proof} 

Now we describe the construction of the action. 
Take $X_{0} = Z$ and $G_{0} = \mathbb{Z}$ with the action $G_{0}\curvearrowright X_{0}$ given by the homeomorphism $\zeta$. 
For each $n\geq 0$, take 
\begin{equation} 
X_{n+1} = X_{n}^{\times \ensuremath{\mathbb{Z}}} \times Z_{n+1}, 
\end{equation} 
where $Z_{n+1}$ is a distinguished copy of $Z$. 
Assuming that $G_{n} \curvearrowright X_{n}$ has been defined, we define recursively an action of $G_{n+1} = \mathfrak{C}(G_{n}) = \mathfrak{D}(G_{n})\rtimes_{\sigma_{G_{n}}}\ensuremath{\mathbb{Z}}$ on $X_{n+1}$. 
For every $((x_{k})_{k}, z)\in X_{n+1}$ and group element $g = (g_{k})_{k}\in \mathfrak{D}(G_{n})$, we define 
\begin{equation} 
\begin{split} 
g\cdot ((x_k)_{k}, z) &=((g_k\cdot x_k)_{k}, z), \\
\quad t_{n+1}\cdot ((x_k)_{k}, z) &=((x_{k+1})_{k}, \zeta_{n+1}(z)), 
\end{split} 
\end{equation} 
where $t_{n+1}$ is the generator of the acting copy of $\mathbb{Z}$ in $G_{n+1}$, and $\zeta_{n+1}$ denotes the homeomorphism on $Z_{n+1}$ given by Proposition~\ref{space_Z}. 
This defines an action $G_{n+1}\curvearrowright X_{n+1}$, as the covariance condition is satisfied: 
\begin{equation}  
   t_{n+1} g t_{n+1}^{-1} \cdot ((x_k)_{k}, z)=((g_{k+1}\cdot x_k)_{k},z) = \sigma_{G_{n}}(g)\cdot  ((x_k)_{k}, z).
\end{equation}

Define the bonding map \(\pi_n\colon X_{n+1}\to X_n\) by \(\pi_n((x_k)_{k},z)=x_0\), and set 
\begin{equation} 
X=\varprojlim\bigl(X_0\xleftarrow{\pi_0}X_1\xleftarrow{\pi_1}X_2\xleftarrow{\pi_2}\cdots\bigr).
\end{equation} 
Note that \(\pi_n\) is \emph{equivariant} with respect to the embedding $\delta_{G_{n}}$. 
Indeed, for any $g \in G_{n}$, since the element \(\delta_{G_n}(g)\) acts by \(g\) on the \(0\)-th coordinate of $X_{n+1}$ and trivially elsewhere, we have 
\begin{equation} 
\pi_n(\delta_{G_n}(g)\cdot ((x_k)_k, z))=g\cdot x_0=g\cdot\pi_n((x_k)_k,z), 
\end{equation} 
for any $((x_k)_{k},z)\in X_{n+1}$. 
Our construction induces an action $G\curvearrowright X$, where 
\begin{equation} 
G = \varinjlim (G_{n}, \delta_{G_{n}}) = \mathfrak{B}(\mathbb{Z})
\end{equation} 
is countable, amenable, torsion-free and acyclic by Proposition~\ref{prop:B(G)-is-nice} and Corollary~\ref{cor:B-acyclic}. 

By \cite[Lemma 2.5.2]{CEL_2015}, $\pi_{n}$ and $\delta_{G_{n}}$ induce unital, injective $^*$-homomorphisms \mbox{$\varphi_{n} \colon C(X_n)\rtimes_{r} G_n \rightarrow C(X_{n+1})\rtimes_{r} G_{n+1}$}. A standard density argument  implies that 
\begin{equation} 
\label{ind_lim_struc} 
C(X)\rtimes_{r} G\cong\varinjlim (C(X_n)\rtimes_{r} G_n, \varphi_{n}). 
\end{equation} 
In the rest of the section, we show that $C(X) \rtimes_{r} G$ is classifiable (the non-trivial parts being simplicity and $\mathcal{Z}$-stability), admits a unique trace and has the same $K$-theory as the complex numbers. 

\begin{thm}
\label{simple_cp} 
The crossed product $C(X)\rtimes_{r} G$ is simple. 
\end{thm} 

\begin{proof} 
We show simplicity of $C(X_n)\rtimes_{r} G_n$, so that simplicity of $C(X)\rtimes_{r} \ensuremath{\mathfrak{B}}(\mathbb{Z})$ follows by taking inductive limits. It thus suffices to show that $G_{n}\curvearrowright X_{n}$ are free and minimal for all $n\geq 0$ (Proposition~\ref{action_prop}~\ref{simple}).  The base case $n=0$ is clear as the action $\zeta\colon \mathbb{Z} \curvearrowright Z$ is free and minimal by Proposition~\ref{space_Z}. Suppose $G_{n}\curvearrowright X_{n}$ is free and minimal for some $n\geq 0$. 

We first prove that $G_{n+1}\curvearrowright X_{n+1}$ is free. Let $((g_{k})_{k}, t_{n+1}^{m}) \in G_{n+1}$, and suppose it fixes some point $((x_{k})_{k}, z)\in X_{n+1}$. Looking at the $Z_{n+1}$-coordinate, we obtain $\zeta^{m}(z) = z$. By freeness of $\zeta$, it follows that $m=0$. 
Thus $((g_{k})_{k}, 0_{\mathbb{Z}})$ fixes $((x_{k})_{k}, z)$, and hence $g_{k}\cdot x_{k}=x_{k}$ for every $k\in \mathbb{Z}$. By freeness of $G_{n}\curvearrowright X_{n}$, we have $g_{k} = 1$ for every $k\in \mathbb{Z}$. 

Now we show minimality of $G_{n+1}\curvearrowright X_{n+1}$. Given any point $((x_{k})_{k}, z)\in X_{n+1}$, it suffices to show that its orbit intersects non-trivially with any basic open set of the form 
\begin{equation}
   \ensuremath{\mathcal{U}}=
   \Bigl(\prod_{k\in F}U_k\Bigr)
   \times
   \Bigl(\prod_{k\notin F}X_n\Bigr)
   \times V, 
\end{equation} 
for a finite subset $F\subseteq \mathbb{Z}$, non-empty open subsets $U_{k}\subseteq X_{n}$ for $k\in F$ and a non-empty open subset $V\subseteq Z_{n+1}$. Since $\zeta$ is minimal, there exists \(m\in\ensuremath{\mathbb{Z}}\) with $\zeta ^{m}(z) \in V$. For each $k\in F$, by minimality of $G_{n}\curvearrowright X_{n}$, there exists $g_{k}\in G_{n}$ such that $g_k\cdot x_{k+m}\in U_k$. Take $g_{k} = 1_{G_{n}}$ for $k\notin F$. Then 
\begin{equation} 
((g_{k})_{k}, t_{n+1}^{m}) \cdot ((x_{k})_{k}, z) = ((g_{k} \cdot x_{k+m})_k, \zeta ^{m}(z)) \in \mathcal{U}. 
\end{equation} 
Thus the orbit of $((x_{k})_{k}, z)$ is dense. 
\end{proof} 

Next we prove $\mathcal{Z}$-stability of the crossed product through the central sequence characterization (see Corollary~\ref{Z_commute_lem}). 

\begin{thm}
\label{Z_stable_trace_cp} 
The crossed product $C(X)\rtimes_{r} G$ is $\ensuremath{\mathcal{Z}}$-stable. 
\end{thm} 

\begin{proof} 
Note first that $C(X_0)\rtimes_{r} G_0 = C(Z)\rtimes_{r} \mathbb{Z}$ is $\ensuremath{\mathcal{Z}}$-stable, as it is induced by a free, minimal action of $\ensuremath{\mathbb{Z}}$ on a compact metric space with finite covering dimension (see \cite[Theorem~B]{tomswinter}).\footnote{This observation was also made in the work of Deeley, Putnam and Strung; see \cite[Corollary 3.3]{DPS-jiang-su}.}
Since the connecting maps $\varphi_{n}$ are unital, injective $^*$-homomorphisms, $C(X_n)\rtimes_{r} G_n$ contains a unital copy of $\ensuremath{\mathcal{Z}}$ for every $n\geq 0$. 

For every $n \geq 0$, let $v_{n+1}$ be the unitary in $C(X_{n+1})\rtimes_{r} G_{n+1}$ corresponding to $t_{n+1}$ in $G_{n+1}$. We claim that $\varphi_n(C(X_n)\rtimes_{r} G_n)$ and $\operatorname{Ad}_{v_{n+1}}\circ \varphi_n(C(X_n)\rtimes_{r} G_n)$ commute inside $C(X_{n+1})\rtimes_{r} G_{n+1}$. It suffices to verify the following four relations: 
\begin{equation} 
\begin{split} 
    & [x,\operatorname{Ad}_{v_{n+1}} (y)]=0,\qquad [u_{g}, \operatorname{Ad}_{v_{n+1}} (u_{h})] = 0, \\ 
    & [x, \operatorname{Ad}_{v_{n+1}} (u_{h})]=0,\qquad [u_{g}, \operatorname{Ad}_{v_{n+1}} (y)]=0, 
\end{split} 
\end{equation} 
for every $x,y\in \varphi_{n}(C(X_n)) \subseteq C(X_{n+1})$ and every $g,h\in \delta_{G_n}(G_n)\leq G_{n+1}$. The first is immediate from the fact that $v_{n+1}$ normalizes $C(X_{n+1})$. The second follows from the fact that $\delta_{G_n}(G_n)$ and $\operatorname{Ad}_{t_{n+1}}\circ \delta_{G_n}(G_n)$ are contained in different summands inside $\mathfrak{D}(G_n)$ and thus commute. 
The third follows from the observation that $\operatorname{Ad}_{t_{n+1}}\circ \delta_{G_n}(G_n)$ acts trivially in the $0$-th coordinate of $X_{n+1}$ and the fourth assertion follows similarly. 

For every $n\geq 0$, $\operatorname{Ad}_{v_{n+1}}\circ \varphi_n(C(X_n)\rtimes_{r} G_n)$ contains a unital copy of $\ensuremath{\mathcal{Z}}$ as this is the case for $C(X_n)\rtimes_{r} G_n$. By our claim above, it follows that for every $n\geq 0$, there is a unital embedding $\iota_{n}\colon \ensuremath{\mathcal{Z}}\rightarrow C(X_{n+1})\rtimes_{r} G_{n+1}$ commuting with $\varphi_n(C(X_n)\rtimes_{r} G_n)$. By Lemma~\ref{Z_commute_lem}, we conclude that $C(X)\rtimes_{r} G$ is $\ensuremath{\mathcal{Z}}$-stable. 
\end{proof} 

To prove that $C(X)\rtimes_{r} G$ admits a unique trace, we show that $C(X_{n}) \rtimes_{r} G_{n}$ has a unique trace through unique ergodicity of $G_{n} \curvearrowright X_{n}$. We start with the following standard lemma. 

\begin{lemma}\label{lem:product-unique}
Let \(\Gamma_i\curvearrowright Y_i\), \(i=1,\dots,r\), be uniquely ergodic actions of countable groups on compact metrizable spaces, with invariant measures \(\mu_i\). Then the product action \(\Gamma_1\times\cdots\times\Gamma_r\curvearrowright Y_1\times\cdots\times Y_r\) is uniquely ergodic, with invariant measure \(\mu_1\times\cdots\times\mu_r\).
\end{lemma}

\begin{proof}
By induction, it is enough to prove the case \(r=2\), which follows from \cite[Lemma 3.1]{David_Li_2024}.
\end{proof} 

The next lemma is a uniquely ergodic formulation of the fact that an ergodic $\Gamma$-system and a trivial $\Gamma$-system only admit the product joining. Here a \emph{joining} is an invariant probability measure on the product of spaces with prescribed marginals. For a probability measure $\eta$ on $Y\times Z$, the \emph{$Z$-marginal} $\lambda$ of $\eta$ is defined by $\lambda (U) = \eta (Y \times U)$ for every Borel set $U\subseteq Z$. 

\begin{lemma}\label{lem:relative-unique} 
Let \(\Gamma\curvearrowright Y\) be a uniquely ergodic action of a countable group $\Gamma$ on a compact metrizable space  $Y$, with invariant measure \(\mu\). Let \(Z\) be compact metrizable. If \(\eta\) is a probability measure on \(Y\times Z\) which is invariant under the action of \(\Gamma\) on the first coordinate, then
\begin{equation}
   \eta=\mu\times\lambda,
\end{equation}
where \(\lambda\) is the \(Z\)-marginal of \(\eta\).
\end{lemma} 

\begin{proof} 
Since $\Gamma \curvearrowright Y$ is uniquely ergodic, it follows that the $Y$-marginal of $\eta$ is $\mu$. Then the lemma follows from \cite[Lemma 6.24]{Glasner_ergodicity_joining} by ergodicity of $\Gamma \curvearrowright Y$. 
\end{proof} 

Combining the lemmas above gives the following theorem. 

\begin{thm}
\label{unique_trace_cp} 
The crossed product $C(X)\rtimes_{r} G$ has a unique trace. 
\end{thm} 

\begin{proof} 
It suffices to show that \(G_n\curvearrowright X_n\) is uniquely ergodic for every $n \geq 0$. Indeed, by unique ergodicity and freeness of finite-stage actions (the latter was shown in the proof of Theorem~\ref{simple_cp}), $C(X_{n}) \rtimes_{r} G_{n}$ has a unique trace $\tau_{n}$ for every $n \geq 0$ (Proposition~\ref{action_prop}~\ref{unique_trace}). Since the connecting maps $\varphi_{n}$ are unital $^*$-homomorphisms, $\tau_{n+1}\circ \varphi_{n}$ is a trace on $C(X_{n}) \rtimes_{r} G_{n}$, and thus equal to $\tau_{n}$ by uniqueness. Therefore $(\tau_{n})_{n}$ induces a trace $\tau$ on $C(X) \rtimes_{r} \ensuremath{\mathfrak{B}}(\ensuremath{\mathbb{Z}})$. Moreover, every trace on the inductive limit restricts to $\tau_{n}$ at each stage, so density of the union of images from each stage implies that $\tau$ is also the unique trace.

Note that the base case action $G_{0}\curvearrowright X_{0}$ given by $\zeta$ is uniquely ergodic by Proposition~\ref{space_Z}. We denote by $\mu_{0}$ the invariant measure for the action. Suppose \(G_n\curvearrowright X_n\) is uniquely ergodic, with the invariant measure \(\mu_n\) for some $n \geq 0$. Let \(\nu_{n+1}\) be the unique \(\zeta_{n+1}\)-invariant probability measure on \(Z_{n+1}\), and put
\begin{equation}
\mu_{n+1}=\mu_n^\ensuremath{\mathbb{Z}}\times\nu_{n+1}.
\end{equation}
This measure is invariant; we show that it is the only invariant probability measure. 

Let \(\eta\) be a \(G_{n+1}\)-invariant probability measure on \(X_{n+1}=X_n^\ensuremath{\mathbb{Z}}\times Z_{n+1}\), and let \(\lambda\) be its \(Z_{n+1}\)-marginal. Fix a non-empty finite set \(F\subset\ensuremath{\mathbb{Z}}\), and let \(\eta_F\) be the marginal of \(\eta\) on \(X_n^F\times Z_{n+1}\). The subgroup of \(\mathfrak{D}(G_n)\) supported on \(F\) is \(G_n^F\), acting coordinatewise on \(X_n^F\) and trivially on \(Z_{n+1}\). Thus \(\eta_F\) is invariant under \(G_n^F\) on the first factor. By Lemma~\ref{lem:product-unique}, the action \(G_n^F\curvearrowright X_n^F\) is uniquely ergodic, with invariant measure \(\mu_n^F\). Lemma~\ref{lem:relative-unique} therefore gives
\begin{equation}
   \eta_F=\mu_n^F\times\lambda .
\end{equation}
Since this holds for every finite \(F\subset\ensuremath{\mathbb{Z}}\), the finite-coordinate marginals of \(\eta\) are those of \(\mu_n^\ensuremath{\mathbb{Z}}\times\lambda\), so
\begin{equation}
   \eta=\mu_n^\ensuremath{\mathbb{Z}}\times\lambda.
\end{equation}
Finally, invariance under \(t_{n+1}\) gives \(\zeta_{n+1*}\lambda=\lambda\). By unique ergodicity of \(\zeta_{n+1}\), \(\lambda=\nu_{n+1}\). Hence \(\eta=\mu_{n+1}\). 
\end{proof}

The last ingredient is to show that $X$ is a point-like space.

\begin{prop}\label{prop:X-properties}
The space \(X\) is compact, metrizable and the inclusion of the unit $\iota\colon\ensuremath{\mathbb{C}}\hookrightarrow C(X)$ induces an isomorphism in $K$-theory.
\end{prop}

\begin{proof}
    As a topological space, $X$ is homeomorphic to the countably infinite product
    \begin{equation}
        X\cong Z^{\times \ensuremath{\mathbb{Z}}}
    \end{equation}
    and is therefore a compact metric space. 
    By \cite[Corollary~1.12]{DPS-jiang-su}, $C(Z)$ is $KK$-equivalent to $\ensuremath{\mathbb{C}}$, and 
    so the inclusion of the unit $\iota_Z\colon \ensuremath{\mathbb{C}}\hookrightarrow C(Z)$ induces an isomorphism in $K$-theory since $C(Z)$ admits a character. 
    
    It follows from the K\"unneth formula that the same is true for any finite product $Z^{\times n}$.
    Finally, $C(X)$ is isomorphic to the inductive limit 
    \begin{equation}
        C(X)\cong \varinjlim C(Z^n)
    \end{equation}
    with unital connecting maps.
    By continuity of $K$-theory, the inclusion of the unit $\iota_X\colon\ensuremath{\mathbb{C}}\hookrightarrow C(X)$ induces an isomorphism in $K$-theory.    
\end{proof}

We now give the proof for Theorem~\ref{introthm:Z_realize}.

\begin{proof}[Proof of Theorem~\ref{introthm:Z_realize}] 
By Theorem~\ref{simple_cp}, Theorem~\ref{unique_trace_cp} and Theorem~\ref{Z_stable_trace_cp}, it follows that $C(X)\rtimes_{r} G$ is simple, $\ensuremath{\mathcal{Z}}$-stable and has a unique trace. Moreover, the crossed product is separable, unital and nuclear by Proposition~\ref{action_prop}. 
The fact that $C(X)\rtimes_{r} G$ satisfies the UCT follows from \cite[Proposition~10.7]{tu}. By Corollary~\ref{cor:K-theory}, the $K$-theory of $C(X)\rtimes_{r} G$ is identical to that of $\ensuremath{\mathcal{Z}}$ and since the class of the unit generates $K_0(C(X)\rtimes_{r} G)$, the pairing maps must also agree.

In conclusion, $C(X)\rtimes_{r} G$ is a unital, simple, separable, nuclear and $\ensuremath{\mathcal{Z}}$-stable $C^*$-algebra satisfying the UCT, and having the same $K$-theory, trace space and pairing map as $\ensuremath{\mathcal{Z}}$. By the unital classification theorem (cf. \cite[Theorem A]{classification}), it follows that $C(X)\rtimes_{r} G\cong \ensuremath{\mathcal{Z}}$.
\end{proof}

\begin{rmk}
    It is possible to produce a crossed product model of \ensuremath{\mathcal{Z}} without invoking the full force of Theorem~\ref{Z_stable_trace_cp}. 
    From Theorems~\ref{simple_cp}, \ref{unique_trace_cp} and Corollary~\ref{cor:K-theory}, the crossed product $C(X)\rtimes_r G$ is unital, simple, separable, nuclear, satisfies the UCT and has the same invariant as the complex numbers. 
    Moreover, it contains a unital copy of \ensuremath{\mathcal{Z}} as it admits a unital embedding of $C(X_{0}) \rtimes_{r} G_{0}$, which is $\ensuremath{\mathcal{Z}}$-stable. 
    Thus, its infinite tensor power has the same properties and is also \ensuremath{\mathcal{Z}}-stable.
    It is evidently still a crossed product, thus establishing Theorem~\ref{introthm:Z_realize}.
   
    We would like to thank W. Winter for pointing out this argument to us.
\end{rmk}

\section{\texorpdfstring{$\mathcal{O}_\infty$}{O-infinity} as a crossed product}
\label{sec:oinfty} 

In this section, we realize $\mathcal{O}_\infty$ as a crossed product by constructing topologically amenable actions of acyclic, word-hyperbolic groups on the Hilbert cube (which is point-like). The resulting crossed product is thus traceless and has the same $K$-theory as the complex numbers. 
The main dynamical input is the following theorem, which is first stated for hyperbolic groups, and will later be applied to an acyclic hyperbolic group.

\begin{thm}
\label{thm:nonamenable-acyclic-point-like-action}
Every non-elementary, torsion-free, word-hyperbolic group admits a minimal,
topologically free, topologically amenable action on the Hilbert cube.
\end{thm}

Fix a group $H$ as in the theorem, and let $X=\partial H$ be its Gromov
boundary. The space $X$ is infinite, compact, metrizable, and has no
isolated points \cite[Theorem~2.28]{KapovichBenakli2002}. We write the
boundary action as $(s,x)\mapsto sx$. Since $H$ has no nontrivial finite
normal subgroup, this action is minimal, strongly proximal, and
topologically free \cite[Proposition~4.1]{AbbottDahmani2019}. It is also
topologically amenable
\cite[Theorem~3.15 and Remark~3.16]{Kaimanovich2004}.
Recall that an action \(H\curvearrowright X\) on a compact, Hausdorff space is
\emph{strongly proximal} if, for every measure \(\mu\in P(X)\), the closure of the orbit \(H\mu\) contains a point mass.

$X$ is, however, not point-like. Therefore, we look at the action on $Q=P(X)$, the space of Borel probability measures on $X$ equipped with the weak-$*$ topology. 
The space $Q$ is infinite-dimensional, compact, convex, and metrizable. Keller's theorem
\cite{Keller1931} therefore identifies $Q$ with the Hilbert cube. 
In particular, $Q$ is contractible and serves as the desired point-like space. 
Choose a sequence $(f_n)_{n}$ that is dense in the unit ball of $C(X)$, and equip $Q$ with the metric 
\begin{equation} 
\label{eq:oinfty-metric}
 d(\mu,\nu)=\frac12\sum_{n=1}^{\infty}2^{-n}
                |\mu(f_n)-\nu(f_n)|.
\end{equation} 
This metric induces the weak-$*$ topology on $Q$, and its diameter is at most $1$.

Let $\alpha\colon H\curvearrowright Q$ be the push-forward action.
For $\mu\in Q$ and $f\in C(X)$, write $\mu(f)=\int_X f\,d\mu$, and the action is given by 
\begin{equation}
 (\alpha_s\mu)(f)=\int_X f(sx)\,d\mu(x).
\end{equation} 
We record several properties of $\alpha$ that will be used throughout. A subset of $Q$ is called \emph{syndetic} if finitely many of its $\alpha$-translates cover $Q$. 

\begin{lemma}
\label{lemma:properties-of-alpha}
The action $\alpha$ is topologically free and topologically amenable.
For every non-empty open set $O\subseteq X$ and every $0<\eta<1$, the set
\begin{equation}
 W_{O,\eta}=\{\mu\in Q \ | \ \mu(O)>\eta\}
\end{equation}
is open and syndetic.
\end{lemma}

\begin{proof}
Topological amenability passes from $X$ to $P(X)$ by
\cite[Lemma~4.4]{KalantarKennedy2017}.
To prove topological freeness, fix $s\ne 1$ and a measure $\mu$ fixed by $\alpha_s$. 
By topological freeness of the action $H \curvearrowright X$, choose $x\in X$ with $sx\ne x$. The measures
\mbox{$\mu_t=(1-t)\mu+t\delta_x$}, for $0<t<1$, converge to $\mu$ as $t\to0$,
and
\begin{equation}
 \alpha_s\mu_t-\mu_t=t(\delta_{sx}-\delta_x)\ne0.
\end{equation}
Thus the fixed-point set of $\alpha_s$ has empty interior.

The set $W_{O, \eta}$ is open as $\mu\mapsto \mu(O)$ is lower-semicontinuous. 
We prove that it is syndetic. 
For every \(\mu\in Q\), strong proximality gives
a point mass \(\delta_x\) in the orbit closure of \(\mu\), and minimality
moves \(x\) into \(O\). 
Hence
\(\bigcup_{g\in H}\alpha_g(W_{O,\eta})=Q\). Since \(W_{O,\eta}\) is open,
compactness gives a finite subcover. 
\end{proof} 

The action $\alpha$ has all the properties required in Theorem~\ref{thm:nonamenable-acyclic-point-like-action} except minimality: the space of point masses forms a proper closed invariant subspace. 
To overcome this, we follow the approach of Glasner--Weiss \cite{glasner1979construction}, and find a minimal action with the rest of the desired properties in the closure of the conjugacy class of $\alpha$, 
\begin{equation} 
 \mathcal A=
 \overline{\{h\alpha h^{-1} \ | \ h\in\operatorname{Homeo}(Q)\}}
 \subseteq\operatorname{Act}(H,Q). 
\end{equation} 
The topology on $\operatorname{Act}(H,Q)$ is uniform convergence on
finite subsets of $H$, under which it is completely metrizable. 
Consequently, the closed subspace $\mathcal A$ is a Baire space; 
see also \cite[pp.~76--77]{suzuki2017construction}. 
We describe briefly the proof strategy. 

For every non-empty open subset $U\subseteq Q$, define 
\begin{equation}
 \mathcal E(U)=
 \left\{\beta\in\mathcal A \ \Big| \
           Q=\bigcup_{g\in H}\beta_g(U)\right\}.
\end{equation} 
The sets $\mathcal E(U)$ are open; see the beginning of the proof of Theorem~\ref{thm:nonamenable-acyclic-point-like-action}. 
More importantly, membership in $\mathcal E(U)$ means that every orbit meets $U$.
Thus an action is minimal if and only if it is in the $\mathcal E(U)$ for every non-empty open subset. 
Our main task is therefore to prove $\mathcal{E}(U)$ is dense in $\mathcal{A}$ for every non-empty open subset $U\subseteq Q$. The Baire category theorem will then give a dense $G_{\delta}$ set of minimal actions in $\mathcal{A}$. 

It suffices to show that $\alpha\in\overline{\mathcal E(U)}$ for every non-empty open $U$. 
Indeed, assume that this is the case and fix $k\in\operatorname{Homeo}(Q)$. 
Conjugation by $k$ carries $\mathcal E(k^{-1}U)$ onto $\mathcal E(U)$,
so the asserted approximation of $\alpha$ also gives an approximation of
$k\alpha k^{-1}$. 
Therefore $k\alpha k^{-1} \in \overline{\mathcal{E} (U)}$, and density of such conjugates in $\mathcal A$ implies that every $\mathcal{E}(U)$ is dense in $\mathcal{A}$. 
The following perturbation lemma proves exactly that $\alpha \in \overline{\mathcal E(U)}$ for any given non-empty open $U$.

\begin{lemma} 
\label{lem:nonamenable-acyclic-point-like-perturbation} 
Let $U\subseteq Q$ be non-empty and open. For every finite subset $E\subseteq H$ 
and every $\varepsilon>0$, there is some $h\in\operatorname{Homeo}(Q)$ such that every $h\alpha h^{-1}$-orbit meets $U$ and 
\begin{equation} 
 \sup_{\mu\in Q}
 d(h\alpha_s h^{-1}\mu,\alpha_s\mu)<\varepsilon,
 \qquad s\in E . 
\end{equation} 
\end{lemma} 

To construct $h$, we use the fact that every $\alpha$-orbit meets each $W_{O,\eta}$. 
It therefore suffices to arrange $h(W_{O,\eta})\subseteq U$ for some open subset $O\subseteq X$ and $\eta>0$. 
To keep the conjugated action close to $\alpha$, we also want $h$ to
almost commute with the finitely many maps $\alpha_s$, $s\in E$.
We first meet these requirements with a continuous affine surjection
$p\colon Q\to Q$. 
Its fibres will be compact convex sets, allowing us to approximate $p$ by a homeomorphism using the Hilbert cube approximation theorem. 

We construct $p$ by first choosing a continuous surjection $q\colon X\rightarrow Q$ and then setting $p\mu = \int_{X} q(x) d\mu(x)$. Thus $q$ prescribes the values of $p$ on the point masses, which are the extreme points in $Q$. 
Approximate equivariance and surjectivity pass from $q$ to $p$. 
Moreover, for a fixed $\mu_0\in U$, if $q(x)=\mu_0$ for all $x\in O$ for some open subset $O\subseteq X$, then for any measure $\mu$ whose mass is concentrated in $O$, the averaging $p\mu$ is close to $\mu_{0}$.

\begin{lemma}
\label{lem:nonamenable-acyclic-point-like-surjection}
Let $\mu_0\in Q$, let $E\subseteq H$ be finite, and let $\kappa>0$.
There are a non-empty open set $O\subseteq X$ and a continuous surjection
$q\colon X\to Q$ such that $q(x)=\mu_0$ for all $x\in O$ and
\begin{equation}
 \sup_{x\in X}d(q(sx),\alpha_s q(x))<\kappa,
 \qquad s\in E .
\end{equation}
\end{lemma}

We first use this lemma to prove the perturbation lemma (Lemma~\ref{lem:nonamenable-acyclic-point-like-perturbation}), and only then give its proof, in order to streamline the argument. 

\begin{proof}[Proof of Lemma~\ref{lem:nonamenable-acyclic-point-like-perturbation}]
Choose $\mu_0\in U$ and $0<r<1$ with $B(\mu_0,2r)\subseteq U$.
Apply Lemma~\ref{lem:nonamenable-acyclic-point-like-surjection} with
error $\kappa = \varepsilon/3$, and let $O\subseteq X$ and $q\colon X\rightarrow Q$ be the resulting open set and surjective map. 

We define $p\colon Q\to Q$ by averaging the values of $q$:
\begin{equation}
\label{eq:oinfty-average}
 (p\mu)(f)=\int_X q(x)(f)\,d\mu(x),
 \qquad f\in C(X) .
\end{equation}
Note that $p\mu$ is indeed a probability measure, since its 
right-hand side is positive, linear in $f$, and takes the value $1$ on the
constant function $1$. It also shows that $p$ is continuous and affine.
Moreover, $p(\delta_x)=q(x)$, so $p$ is onto. Identifying $x$ with $\delta_x$, the map $p$ is the continuous
affine extension of $q$ to $Q$. 

We will use the following estimate for averaging twice. 
For continuous maps $a,b\colon X\to Q$, integrals being understood as
in~\eqref{eq:oinfty-average}, the definition of $d$ in~\eqref{eq:oinfty-metric} gives
\begin{equation}
 \label{eq:oinfty-integral-estimate}
\begin{aligned}
 d\left(\int_X a\,d\mu,\int_X b\,d\mu\right)
 &=\frac12\sum_{n=1}^{\infty}2^{-n}
       \left|\int_X\bigl(a(x)(f_n)-b(x)(f_n)\bigr)
                      \,d\mu(x)\right| \\
 &\leq\int_X d(a(x),b(x))\,d\mu(x), 
\end{aligned}
\end{equation} 
for every $\mu \in Q$. 
In the last step, we apply the triangle inequality, and interchange a
non-negative sum and an integral.

We show approximate equivariance of $p$. The push-forward formula and affinity of $\alpha_s$ give
\begin{equation}
 p(\alpha_s\mu)=\int_X q(sx)\,d\mu(x),
 \quad
 \alpha_s(p\mu)=\int_X\alpha_s q(x)\,d\mu(x), \qquad s\in E, \mu \in Q. 
\end{equation}
Therefore an application of~\eqref{eq:oinfty-integral-estimate} yields
\begin{equation}
\label{eq:oinfty-p-equivariance}
 d(p\alpha_s\mu,\alpha_s p\mu)
 \leq\int_X d(q(sx),\alpha_s q(x))\,d\mu(x)
 <\varepsilon/3,
 \qquad s\in E, \mu \in Q.
\end{equation} 

Take $W=W_{O,1-r}$ consisting of measures assigning almost all mass to $O$, and $W$ is syndetic by
Lemma~\ref{lemma:properties-of-alpha}. We prove that $p(W) \subseteq B(\mu_{0}, r)$, leaving room to approximate $p$ by a homeomorphism while keeping the image of $W$ inside $U$. 
For $\mu\in W$, apply \eqref{eq:oinfty-integral-estimate} with $a=q$ and $b$ the constant map with value $\mu_0$. Since $q(x)=\mu_0$ for $x\in O$ and
$\operatorname{diam}(Q)\leq1$, we obtain
\begin{equation}
\label{eq:oinfty-p-image}
 d(p\mu,\mu_0)\leq\int_Xd(q(x),\mu_0)\,d\mu(x)
             \leq1-\mu(O)<r.
\end{equation}

We now use the Hilbert cube structure of $Q$ through the following consequence
of Chapman's approximation theorem: a continuous
surjection between Hilbert cubes with contractible fibres is uniformly
approximable by homeomorphisms \cite{Chapman1976}. 
Since $p$ is a continuous affine surjection, it follows that each fibre $p^{-1}(\nu)$, $\nu \in Q$, is non-empty, compact, and convex, hence contractible for every $\nu \in Q$. 
The theorem therefore applies to $p$.

By uniform continuity of the finitely many maps $\alpha_s$, choose a
homeomorphism $h$ uniformly close enough to $p$ so that
\begin{equation}
\label{eq:oinfty-h-close-to-p}
 \sup_{\nu\in Q}d(h\nu,p\nu)<\min\{r,\varepsilon/3\},
 \qquad
 \sup_{\nu\in Q}d(\alpha_s h\nu,\alpha_s p\nu)<\varepsilon/3,
 \quad s\in E .
\end{equation}
For $s\in E$ and $\mu\in Q$, the three errors can now be bounded:
\begin{equation}
\begin{aligned}
 d(h\alpha_s\mu,\alpha_s h\mu)
 &\leq d(h\alpha_s\mu,p\alpha_s\mu)
       +d(p\alpha_s\mu,\alpha_s p\mu)\\
 &\hspace{1.3em}+d(\alpha_s p\mu,\alpha_s h\mu)
 <\varepsilon.
\end{aligned}
\end{equation}
Substituting $\mu=h^{-1}\nu$ gives the required estimate for
$h\alpha_s h^{-1}$.

Since $p(W)\subseteq B(\mu_{0}, r)$, estimate~\eqref{eq:oinfty-h-close-to-p} also gives 
$h(W)\subseteq B(\mu_0,2r)\subseteq U$.
Given $\nu\in Q$, choose $g\in H$ with
$\alpha_g(h^{-1}\nu)\in W$, using Lemma~\ref{lemma:properties-of-alpha}.
Then
\begin{equation}
 (h\alpha_g h^{-1})(\nu)\in h(W)\subseteq U.
\end{equation}
Thus every orbit of the conjugated action meets $U$.
\end{proof}

It remains to construct $q$. Our use of amenability to construct approximately equivariant averages is inspired by Suzuki's extension of this approach to amenable actions \cite{suzuki2017construction}. 

\begin{proof}[Proof of Lemma~\ref{lem:nonamenable-acyclic-point-like-surjection}]

By amenability of the action $H\curvearrowright X$, there is a continuous map
$m\colon X\to\operatorname{Prob}_{f}(H)$, whose supports are all
contained in one finite set $F\subseteq H$, such that
\begin{equation}
\label{eq:oinfty-amenability}
 \sup_{x\in X}\|m(sx)-s m(x)\|_1<\kappa,
 \qquad s\in E .
\end{equation}
For $g\in H$, write $m_g(x)=m(x)(g)$. In particular, $m_g=0$ if
$g\notin F$.

Given a continuous map $\theta\colon X\to Q$, define 
\begin{equation}
\label{eq:oinfty-q}
 q(x)=\sum_{g\in F}m_g(x)\,
                     \alpha_g\theta(g^{-1}x).
\end{equation}
For any map $\theta$, estimate~\eqref{eq:oinfty-amenability} will
give approximate equivariance of $q$. We choose $\theta$ to ensure that $q$ is surjective and has constant value $\mu_{0}$ on a non-empty set $O$. 

Topological freeness allows us to choose $x_0\in X$ such that the points
$g^{-1}x_0$, $g\in F$, are distinct: we avoid the fixed-point
sets of the finitely many nonidentity elements
$g_1g_2^{-1}$ with $g_1,g_2\in F$.
Next, choose an open neighbourhood $V$ of $x_0$ such that the sets
\begin{equation}
 g^{-1}\overline V,\qquad g\in F,
\end{equation}
are pairwise disjoint. Since $X$ has no isolated points, we can choose
a non-empty open set $O$ and a Cantor set $R$ inside $V$ with
$\overline O\subseteq V$ and $\overline O\cap R=\varnothing$.
Choose a continuous surjection $\rho\colon R\to Q$. Such a map exists
because every non-empty compact metrizable space is a continuous image of
the Cantor set.

Now, we define the value of $\theta$ on $g^{-1}\overline{O} \cup g^{-1}R$, $g\in F$, by
\begin{equation}
\label{eq:oinfty-theta-def}
\begin{aligned}
 \theta(g^{-1}x)&=\alpha_{g^{-1}}\mu_0,
 && x\in\overline O,\ g\in F ,\\
 \theta(g^{-1}x)&=\alpha_{g^{-1}}\rho(x),
 &&x\in R,\ g\in F.
\end{aligned}
\end{equation}
The sets on which these values are defined are pairwise disjoint and
closed, so \eqref{eq:oinfty-theta-def} defines a continuous map from their union.
They also extend to a continuous map $\theta\colon X\to Q$. Indeed, choose a
homeomorphism from $Q$ onto $[0,1]^{\mathbb N}$ and apply Tietze's
extension theorem to each coordinate \cite[Theorem~35.1(a)]{Munkres2000}.

For $x\in\overline O$, every term
$\alpha_g\theta(g^{-1}x)$ in~\eqref{eq:oinfty-q} is $\mu_0$.
For $x\in R$, every such term is $\rho(x)$.
Since the coefficients sum to $1$, it follows that
$q(x)=\mu_0$ for all $x\in\overline O$ and $q|_R=\rho$.
In particular, $q$ is onto.

To check approximate equivariance, fix $s\in E$.
Since $\alpha_s$ is affine, we can apply it term by term
in the formula \eqref{eq:oinfty-q} for $q(x)$ and $q(sx)$.
Reindexing by $h=sg$, we have
$\alpha_s\alpha_g=\alpha_h$ and
$g^{-1}x=h^{-1}sx$.
The coefficients in the two sums vanish outside $F$ and $sF$,
respectively, so we may write both sums over $F\cup sF$:
\begin{equation}
\begin{aligned}
 q(sx)
 &=\sum_{h\in F\cup sF}m_h(sx)\,
                         \alpha_h\theta(h^{-1}sx),\\
 \alpha_s q(x)
 &=\sum_{h\in F\cup sF}m_{s^{-1}h}(x)\,
                         \alpha_h\theta(h^{-1}sx).
\end{aligned}
\end{equation}
These expressions average the same points of $Q$ with only the coefficients
being different. For probability vectors $(a_h)$ and $(b_h)$ and measures
$\nu_h\in Q$, the definition of $d$ and the fact that
$\|f_n\|_\infty\leq1$ give
\begin{equation}
\begin{aligned}
 d\left(\sum_{h\in F\cup sF}  a_h\nu_h,\sum_{h\in F\cup sF} b_h\nu_h\right)
 &=\frac12\sum_{n=1}^{\infty}2^{-n}
         \left|\sum_{h\in F\cup sF}(a_h-b_h)\nu_h(f_n)\right|\\
 &\leq\frac12\sum_{h\in F\cup sF}|a_h-b_h|.
\end{aligned}
\end{equation}
Apply this with
$a_h=m_h(sx)$, $b_h=m_{s^{-1}h}(x)$, and
$\nu_h=\alpha_h\theta(h^{-1}sx)$.
Since $(s m(x))(h)=m_{s^{-1}h}(x)$, we obtain
\begin{equation}
 d(q(sx),\alpha_s q(x))
 \leq\frac12\|m(sx)-s m(x)\|_1<\kappa,
\end{equation}
uniformly in $x$, as required.
\end{proof}

So far, we have proved that each $\mathcal E(U)$ is dense. We now complete
the proof of Theorem~\ref{thm:nonamenable-acyclic-point-like-action}.
In fact, the actions with all three required properties form a dense
$G_\delta$ subset of $\mathcal A$.

\begin{proof}[Completion of the proof of
Theorem~\ref{thm:nonamenable-acyclic-point-like-action}]
We first check that $\mathcal E(U)$ is open.
If every $\beta$-orbit meets $U$, compactness gives a finite open cover
$V_1,\ldots,V_n$ of $Q$ and elements $g_1,\ldots,g_n\in H$ such that
$\beta_{g_i}(\overline V_i)\subseteq U$ for every $i$.
These inclusions are stable under small uniform perturbations of the maps
$\beta_{g_i}$. Thus $\mathcal E(U)$ is open.
It is dense in $\ensuremath{\mathcal{A}} $ by Lemma~\ref{lem:nonamenable-acyclic-point-like-perturbation}
and the conjugation argument above. If $\mathcal B$ is a countable base
of non-empty open subsets of $Q$, the minimal actions in $\ensuremath{\mathcal{A}} $ are exactly
those in $\bigcap_{U\in\mathcal B}\mathcal E(U)$, which is therefore a dense $G_{\delta}$ subset of $\mathcal{A}$. 

Topological freeness is also a $G_\delta$ condition.
For each $s\ne1$ and $U\in\mathcal B$, the condition that some point
of $U$ is moved by $\beta_s$ is open. It holds for every conjugate of
$\alpha$, so each of these open sets is dense in $\mathcal A$.
Their intersection is exactly the set of topologically free actions.

For amenability, fix a finite subset $E\subseteq H$ and $\delta>0$. 
Let $\mathcal M(E,\delta)$ consist of those $\beta\in\mathcal A$ for which
there is a continuous map $m\colon Q\to\operatorname{Prob}_{f}(H)$,
with common finite support, satisfying 
\begin{equation} 
 \max_{s\in E}\sup_{\mu\in Q}
       \|m(\beta_s\mu)-s m(\mu)\|_1<\delta.
\end{equation}
For fixed $m$, this inequality is again stable under small uniform
perturbations of $\beta$, so $\mathcal M(E,\delta)$ is open.
It contains every conjugate of $\alpha$ and hence is dense.
Choose increasing finite sets $E_n$ with union $H$.
The topologically amenable actions are exactly
\begin{equation}
 \bigcap_{n=1}^{\infty}\mathcal M(E_n,1/n).
\end{equation}
The Baire category theorem shows that the intersection of these three
dense $G_\delta$ sets is dense in $\mathcal A$. Choosing an action in
this intersection proves the theorem.
\end{proof}

We now apply the construction to an acyclic group.

\begin{thm}
\label{thm:Oinfty-crossed-product}
Every non-elementary, torsion-free, word-hyperbolic, acyclic group $H$
admits an action on
the Hilbert cube $Q$ such that
\(
 C(Q)\rtimes_r H\cong\mathcal O_\infty.
\)
\end{thm}

\begin{proof}
    Let $\beta\colon H\curvearrowright Q$ be an action given by
Theorem~\ref{thm:nonamenable-acyclic-point-like-action}, and set
$A=C(Q)\rtimes_{r,\beta}H$.
For a non-elementary hyperbolic group, a minimal, topologically free,
topologically amenable action on a compact metrizable space has a
crossed product which is a Kirchberg algebra satisfying the UCT
\cite[Theorem~C and Corollary~3.8]{GardellaGeffenKranzNaryshkin2023}.
Thus $A$ is a unital UCT Kirchberg algebra.

Since $Q$ is contractible, the  inclusion
$\mathbb C\to C(Q)$ induces an isomorphism in $K$-theory.
The group $H$ satisfies the Baum--Connes conjecture with coefficients
\cite{lafforgue12}. Together with acyclicity of $H$, Theorem~\ref{thm:KKG} and
Lemma~\ref{lem:baum-connes} give
\begin{equation}
 K_0(A)=\mathbb Z\cdot[1_A],\qquad K_1(A)=0,
\end{equation}
with the unit class preserved by the isomorphisms.
By the Kirchberg--Phillips classification theorem
\cite[Theorem~4.2.4]{Phillips2000},
$A\cong\mathcal O_\infty$.
\end{proof}

Crucially, groups satisfying the hypotheses do exist.
The group in \cite[Proposition~3.2]{Leary2013} is nontrivial,
torsion-free, and acyclic, and its finite presentation complex $K$ admits a
negatively curved metric. The lifted metric makes the universal cover $\widetilde K$ a
hyperbolic geodesic space with compact closed balls
\cite[Theorem~II.4.1(2) and p.~399]{BridsonHaefliger1999}.
The action of $H=\pi_1(K)$ by covering transformations is therefore isometric
and proper: for each compact $C\subseteq\widetilde K$, only finitely
many $g\in H$ satisfy $gC\cap C\ne\varnothing$.
It is cocompact because the quotient is the compact space $K$.
The Milnor--\v{S}varc lemma gives a quasi-isometry between a Cayley
graph of $H$ and $\widetilde K$. Hyperbolicity is preserved by
quasi-isometry, so $H$ is word-hyperbolic
\cite[Proposition~I.8.19 and Theorem~III.H.1.9]{BridsonHaefliger1999}. Finally, it is non-elementary since a torsion-free, elementary hyperbolic group is isomorphic to $\mathbb{Z}$, but $\mathbb{Z}$ is not acyclic.

\begin{rmk}
    There is an alternative approach to obtaining a crossed product model for $\ensuremath{\mathcal{O}} _\infty$ from Theorem~\ref{thm:nonamenable-acyclic-point-like-action}. 
    Instead of appealing to \cite{GardellaGeffenKranzNaryshkin2023} to ensure dynamical comparison, one can do the following:
    Let $\beta$ be the action of a torsion-free, non-elementary hyperbolic, acyclic group $H$ on the Hilbert cube $Q$ produced by Theorem~\ref{thm:nonamenable-acyclic-point-like-action} and let $\alpha$ be the action of $\ensuremath{\mathfrak{B}}(\mathbb{Z})$ on the space $X$ from Section~\ref{sec:model-for-Z}.
    Consider the product action $\beta\times\alpha$ of $H\times \ensuremath{\mathfrak{B}}(\ensuremath{\mathbb{Z}})$ on $Q\times X$.
    The crossed product 
    \begin{equation}
        C(Q\times X)\rtimes_{r,\beta\times\alpha}
\bigl(H\times\mathfrak{B}(\ensuremath{\mathbb{Z}})\bigr)
\cong
\bigl(C(Q)\rtimes_{r,\beta}H\bigr)\otimes\mathcal{Z}
    \end{equation}
    is $\ensuremath{\mathcal{Z}}$-stable and therefore purely infinite by \cite{BK_dichotomy}.
    By the Kirchberg-Phillips classification theorem (\cite{KP_classi}), it is isomorphic to $\ensuremath{\mathcal{O}} _\infty$.
\end{rmk}

\section{Application to  
non-commutative dynamics} 
\label{sec:appli_XLVI} 

The goal of this section is to realize every $\ensuremath{\mathcal{Z}}$-stable \ensuremath{C^{*}}-algebra $A$ as a crossed product
\begin{equation}
    A\cong A\rtimes_{r} G
\end{equation}
for any torsion-free, amenable, acyclic group $G$.
This is achieved by realizing the Jiang--Su algebra $\ensuremath{\mathcal{Z}} $ as a crossed product of a Bernoulli shift action on itself and then tensoring it with the identity action on $A$. An action $\alpha$ on a unital $C^*$-algebra $A$ with $T(A) \neq \emptyset$ is \emph{strongly outer} if for every $g \neq 1_{G}$ and trace $\tau\in T(A)$ satisfying $\tau \circ \alpha_{g} = \tau$, the induced action on $\pi_{\tau}(A)''$ is outer. 

\begin{thm} 
\label{thm:outer-action-on-Z}
Let $G$ be a non-trivial, amenable, countable, torsion-free, acyclic group. 
Let $\beta\colon G \curvearrowright \ensuremath{\mathcal{Z}} \cong \otimes_{G}\ensuremath{\mathcal{Z}}$ be the Bernoulli shift action. 
Then 
$\ensuremath{\mathcal{Z}}\cong \ensuremath{\mathcal{Z}}\rtimes_{r,\beta} G$. 
\end{thm}

\begin{proof} 
    Take the associated crossed product of the Bernoulli action $\beta$, 
    \begin{equation}
        A=\left(\bigotimes_{g\in G}\ensuremath{\mathcal{Z}}\right)\rtimes_{r,\beta} G.
    \end{equation} 
    We will show that $A\cong \ensuremath{\mathcal{Z}}$. 
    Since $G$ is amenable, it satisfies the Baum-Connes conjecture with coefficients. Therefore, by Theorem~\ref{thm:KKG}, 
    \begin{equation}
        K_*(A)\cong K_*(\ensuremath{C^{*}_r}(G))\cong (\ensuremath{\mathbb{Z}}\cdot [1]_0,0)
    \end{equation}
    and both isomorphisms preserve the class of the unit in $K_0$. 
    
    The fact that $A$ is nuclear, $\ensuremath{\mathcal{Z}}$-stable, has a unique trace and satisfies the UCT follows from an argument identical to the one in \cite[Example~9.11]{classification}, which we include for the convenience of the reader. 
    Since nuclearity (see \cite[IV.3.5.2]{blackadar}) and the UCT (\cite[Corollary~9.4]{meyernest}, see also \cite[p.\ 40]{99_open_nuclear}) are preserved by crossed products by torsion-free, amenable groups, $A$ is nuclear and satisfies the UCT. 
    The action on the GNS completion $\bigotimes_{g\in G}\ensuremath{\mathcal{R}}$ is by non-trivial permutations of the index set and so is properly outer (\cite[Theorem 5.2.4]{Anantharaman_Popa}). In particular, the action on the GNS completion is outer and the Bernoulli shift $\beta$ is strongly outer as $\otimes_{g\in G} \mathcal Z\cong \mathcal Z$ admits a unique trace. 
    Therefore, $A$ is simple by \cite[Theorem~3.1]{kishimoto}, $\ensuremath{\mathcal{Z}}$-stable by \cite[Theorem~1.1]{sato} and has a unique trace by \cite[Corollary~1.7]{ursu}.
    The unital classification theorem (cf. \cite[Theorem A]{classification}) finishes the proof. 
\end{proof}

\begin{cor}
    Let $G$ be a torsion-free, amenable, acyclic group.
    Let $A$ be a $\ensuremath{\mathcal{Z}}$-stable \ensuremath{C^{*}}-algebra. 
    Then there exists an outer action $\alpha$ of $G$ on $A$ such that 
    \begin{equation} 
        A\rtimes_{r,\alpha} G\cong A.
    \end{equation}
    In particular, every classifiable \ensuremath{C^{*}}-algebra arises as a crossed product of a classifiable \ensuremath{C^{*}}-algebra by an outer action of a non-trivial amenable group.
\end{cor}

\begin{proof}
    Since $A$ is $\ensuremath{\mathcal{Z}}$-stable, $A\cong A\otimes \ensuremath{\mathcal{Z}}$. 
    Consider the action $\ensuremath{\mathrm{id}}\otimes\beta$ of $G$ on $A\otimes\ensuremath{\mathcal{Z}}$, where $\beta$ is the Bernoulli shift action on $\mathcal Z$. 
    Since $\beta$ is strongly outer as shown in the proof of Theorem~\ref{thm:outer-action-on-Z}, it follows that $\ensuremath{\mathrm{id}}\otimes\beta$ is outer\footnote{This follows from \cite[Corollary~1.12]{kallman} applied to the weak-$*$ closure of $\pi\otimes\pi_\tau(A\otimes\ensuremath{\mathcal{Z}} )$, where $\pi$ is the universal representation of $A$ and $\pi_\tau$ the GNS representation associated to $\tau$. Here $\tau$ is the unique trace on $\ensuremath{\mathcal{Z}}$.}.
    By Theorem~\ref{thm:outer-action-on-Z}, 
    \begin{equation}(A\otimes\ensuremath{\mathcal{Z}})\rtimes_{r,\ensuremath{\mathrm{id}}\otimes\beta}G\cong A\otimes(\ensuremath{\mathcal{Z}}\rtimes_{r,\beta} G)\cong A\otimes\ensuremath{\mathcal{Z}}\cong A,
    \end{equation} 
    which completes the proof. 
\end{proof} 

\newcommand{\etalchar}[1]{$^{#1}$}


\begin{thebibliography}{AMCFQ26}

\bibitem[AD02]{Claire_amenable}
Claire Anantharaman-Delaroche.
\newblock Amenability and exactness for dynamical systems and their
  {$C^\ast$}-algebras.
\newblock {\em Trans. Amer. Math. Soc.}, 354(10):4153--4178, 2002.

\bibitem[AD19]{AbbottDahmani2019}
C.~R. Abbott and F.~Dahmani.
\newblock Property {$P_{\mathrm{naive}}$} for acylindrically hyperbolic groups.
\newblock {\em Math. Z.}, 291(1--2):555--568, 2019.

\bibitem[ADR00]{AnantharamanDelarocheRenault2000}
C.~Anantharaman-Delaroche and J.~Renault.
\newblock {\em Amenable Groupoids}, volume~36 of {\em Monographies de
  L'Enseignement Math{\'e}matique}.
\newblock L'Enseignement Math{\'e}matique, Geneva, 2000.

\bibitem[AMCFQ26]{amcfq}
J.~F. Ariza~Mejía, I.~Chifan, and A.~Fernández~Quero.
\newblock Non-amenable ${C}^*$-superrigid groups that are not
  ${W}^*$-superrigid, 2026.
\newblock Preprint, arXiv:2602.05290.

\bibitem[AP]{Anantharaman_Popa}
C.~Anantharaman and S.~Popa.
\newblock An introduction to {$\mathrm{II}_1$} factors.
\newblock Draft monograph.

\bibitem[AS94]{archboldspielberg}
R.~J. Archbold and J.~S. Spielberg.
\newblock Topologically free actions and ideals in discrete {$C^*$}-dynamical
  systems.
\newblock {\em Proc. Edinburgh Math. Soc. (2)}, 37(1):119--124, 1994.

\bibitem[BBS{\etalchar{+}}19]{BBSTWW}
J.~Bosa, N.~P. Brown, Y.~Sato, A.~Tikuisis, S.~White, and W.~Winter.
\newblock Covering dimension of {$ {C}^*$}-algebras and 2-coloured
  classification.
\newblock {\em Mem. Amer. Math. Soc.}, 257(1233):vii+97, 2019.

\bibitem[BCM04]{BCM2004}
A.~J. Berrick, I.~Chatterji, and G.~Mislin.
\newblock From acyclic groups to the {B}ass conjecture for amenable groups.
\newblock {\em Math. Ann.}, 329(4):597--621, 2004.

\bibitem[BDH80]{baumslag1980topology}
G.~Baumslag, E.~Dyer, and A.~Heller.
\newblock The topology of discrete groups.
\newblock {\em J. Pure Appl. Algebra}, 16(1):1--47, 1980.

\bibitem[Ber89]{Berrick1989}
A.~J. Berrick.
\newblock Universal groups, binate groups and acyclicity.
\newblock In {\em Group theory ({S}ingapore, 1987)}, pages 253--266. de
  Gruyter, Berlin, 1989.

\bibitem[BH99]{BridsonHaefliger1999}
M.~R. Bridson and A.~Haefliger.
\newblock {\em Metric Spaces of Non-Positive Curvature}, volume 319 of {\em
  Grundlehren der mathematischen Wissenschaften}.
\newblock Springer-Verlag, Berlin, 1999.

\bibitem[BK04]{BK_dichotomy}
E.~Blanchard and E.~Kirchberg.
\newblock Non-simple purely infinite {$C^*$}-algebras: the {H}ausdorff case.
\newblock {\em J. Funct. Anal.}, 207(2):461--513, 2004.

\bibitem[BKKO17]{BKKO}
E.~Breuillard, M.~Kalantar, M.~Kennedy, and N.~Ozawa.
\newblock {$C^*$}-simplicity and the unique trace property for discrete groups.
\newblock {\em Publ. Math. Inst. Hautes \'Etudes Sci.}, 126:35--71, 2017.

\bibitem[Bla90]{Blackadar_CAR}
B.~Blackadar.
\newblock Symmetries of the {CAR} algebra.
\newblock {\em Ann. of Math. (2)}, 131(3):589--623, 1990.

\bibitem[Bla06]{blackadar}
B.~Blackadar.
\newblock {\em Operator algebras}, volume 122 of {\em Encyclopaedia of
  Mathematical Sciences}.
\newblock Springer-Verlag, Berlin, 2006.
\newblock Theory of $C^*$-algebras and von Neumann algebras, Operator Algebras
  and Non-commutative Geometry, III.

\bibitem[BO08]{brownozawa}
N.~P. Brown and N.~Ozawa.
\newblock {\em {$C^*$}-algebras and finite-dimensional approximations},
  volume~88 of {\em Graduate Studies in Mathematics}.
\newblock American Mathematical Society, Providence, RI, 2008.

\bibitem[CDAD23a]{cdad1}
I.~Chifan, A.~Diaz-Arias, and D.~Drimbe.
\newblock New examples of {$ W^*$} and {$ C^*$}-superrigid groups.
\newblock {\em Adv. Math.}, 412:Paper No. 108797, 57, 2023.

\bibitem[CDAD23b]{cdad2}
I.~Chifan, A.~Diaz-Arias, and D.~Drimbe.
\newblock {$W^*$} and {$C^*$}-superrigidity results for coinduced groups.
\newblock {\em J. Funct. Anal.}, 284(1):Paper No. 109730, 48, 2023.

\bibitem[CEKN24]{CEKN_Bernoulli}
S.~Chakraborty, S.~Echterhoff, J.~Kranz, and S.~Nishikawa.
\newblock {$K$}-theory of noncommutative {Bernoulli} shifts.
\newblock {\em Math. Ann.}, 388:2671--2703, 2024.

\bibitem[CEL15]{CEL_2015}
J.~Cuntz, S.~Echterhoff, and X.~Li.
\newblock On the {$K$}-theory of the ${C}^*$-algebra generated by the left
  regular representation of an {O}re semigroup.
\newblock {\em J. Eur. Math. Soc. (JEMS)}, 17(3):645--687, 2015.

\bibitem[CELY17]{KKG-book}
J.~Cuntz, S.~Echterhoff, X.~Li, and G.~Yu.
\newblock {\em {$K$}-theory for group {$C^*$}-algebras and semigroup
  {$C^*$}-algebras}, volume~47 of {\em Oberwolfach Seminars}.
\newblock Birkh\"auser/Springer, Cham, 2017.

\bibitem[CFaH20]{CFaH_ample_Kirchberg}
L.~O. Clark, J.~Fletcher, and A.~an~Huef.
\newblock All classifiable {K}irchberg algebras are {$C^*$}-algebras of ample
  groupoids.
\newblock {\em Expo. Math.}, 38(4):559--565, 2020.

\bibitem[CGS{\etalchar{+}}]{classification2}
J.~R. Carrión, J.~Gabe, C.~Schafhauser, A.~Tikuisis, and S.~White.
\newblock Classifying $^*$-homomorphisms {II}.
\newblock In preparation.

\bibitem[CGS{\etalchar{+}}23]{classification}
J.~R. Carrión, J.~Gabe, C.~Schafhauser, A.~Tikuisis, and S.~White.
\newblock Classifying $^*$-homomorphisms {I}: Unital simple nuclear
  {$C^*$}-algebras, 2023.
\newblock arXiv:2307.06480v4, To appear in Acta Math.

\bibitem[Cha76]{Chapman1976}
T.~A. Chapman.
\newblock {\em Lectures on Hilbert Cube Manifolds}, volume~28 of {\em CBMS
  Regional Conference Series in Mathematics}.
\newblock American Mathematical Society, Providence, RI, 1976.

\bibitem[CI18]{ci18}
I.~Chifan and A.~Ioana.
\newblock Amalgamated free product rigidity for group von {N}eumann algebras.
\newblock {\em Adv. Math.}, 329:819--850, 2018.

\bibitem[CIOS23]{CIOS}
I.~Chifan, A.~Ioana, D.~Osin, and B.~Sun.
\newblock Wreath-like products of groups and their von {N}eumann algebras {I}:
  {$ W^\ast $}-superrigidity.
\newblock {\em Ann. of Math. (2)}, 198(3):1261--1303, 2023.

\bibitem[CKR{\etalchar{+}}]{ckrtw}
J.~Curda, S.~Knudby, S.~Raum, H.~Thiel, and S.~White.
\newblock On ${C}^*$-superrigidity of finitely generated torsion free virtually
  abelian groups.
\newblock In preparation.

\bibitem[CMO24]{caprace2024simple}
P.-E. Caprace, A.~Minasyan, and D.~Osin.
\newblock Simple $p$-adic {Lie} groups with abelian {Lie} algebras.
\newblock {\em J. Reine Angew. Math.}, 2024(812):229--256, 2024.

\bibitem[Con76]{connes76}
A.~Connes.
\newblock Classification of injective factors. {C}ases {$\rm{II}\sb{1},$}
  {$\rm{II}\sb{\infty },$} {$\rm{III}\sb{\lambda },$} {$\lambda \not=1$}.
\newblock {\em Ann. of Math. (2)}, 104(1):73--115, 1976.

\bibitem[Dix77]{Dixmier}
J.~Dixmier.
\newblock {\em {$C\sp*$}-algebras}, volume Vol. 15 of {\em North-Holland
  Mathematical Library}.
\newblock North-Holland Publishing Co., Amsterdam-New York-Oxford, 1977.
\newblock Translated from the French by Francis Jellett.

\bibitem[DPS18]{DPS-jiang-su}
R.~J. Deeley, I.~F. Putnam, and K.~R. Strung.
\newblock Constructing minimal homeomorphisms on point-like spaces and a
  dynamical presentation of the {J}iang-{S}u algebra.
\newblock {\em J. Reine Angew. Math.}, 742:241--261, 2018.

\bibitem[DR00]{dykemaradulescu}
K.~J. Dykema and F.~R{\u a}dulescu.
\newblock Compressions of free products of von {N}eumann algebras.
\newblock {\em Math. Ann.}, 316(1):61--82, 2000.

\bibitem[DS01]{dykemashlyakhtenko}
K.~J. Dykema and D.~Shlyakhtenko.
\newblock Exactness of {C}untz-{P}imsner {$C^*$}-algebras.
\newblock {\em Proc. Edinb. Math. Soc. (2)}, 44(2):425--444, 2001.

\bibitem[DV73]{DyerVasquez1973}
E.~Dyer and A.~T. Vasquez.
\newblock Some small aspherical spaces.
\newblock {\em J. Aust. Math. Soc.}, 16:332--352, 1973.

\bibitem[Dyk93]{dykema93}
K.~Dykema.
\newblock Free products of hyperfinite von {N}eumann algebras and free
  dimension.
\newblock {\em Duke Math. J.}, 69(1):97--119, 1993.

\bibitem[EK01]{EK.PJM}
G.~A. Elliott and D.~Kucerovsky.
\newblock An abstract {V}oiculescu-{B}rown-{D}ouglas-{F}illmore absorption
  theorem.
\newblock {\em Pacific J. Math.}, 198(2):385--409, 2001.

\bibitem[EKQR06]{EKQR}
S.~Echterhoff, S.~Kaliszewski, J.~Quigg, and I.~Raeburn.
\newblock A categorical approach to imprimitivity theorems for
  {$C^*$}-dynamical systems.
\newblock {\em Mem. Amer. Math. Soc.}, 180(850):viii+169, 2006.

\bibitem[ER18]{eckhardtraum}
C.~Eckhardt and S.~Raum.
\newblock {$ C^*$}-superrigidity of 2-step nilpotent groups.
\newblock {\em Adv. Math.}, 338:175--195, 2018.

\bibitem[EST22]{EST_KMS}
G.~A. Elliott, Y.~Sato, and K.~Thomsen.
\newblock On the bundle of {KMS} state spaces for flows on a
  {$\mathcal{Z}$}-absorbing {$C^*$}-algebra.
\newblock {\em Comm. Math. Phys.}, 393(2):1105--1123, 2022.

\bibitem[FFLM23]{fournier-facio}
F.~Fournier-Facio, C.~L\"oh, and M.~Moraschini.
\newblock Bounded cohomology and binate groups.
\newblock {\em J. Aust. Math. Soc.}, 115(2):204--239, 2023.

\bibitem[Gar23]{GardamComplexUnits}
G.~Gardam.
\newblock Non-trivial units of complex group rings, 2023.
\newblock Preprint, arXiv:2312.05240.

\bibitem[GGKN23]{GardellaGeffenKranzNaryshkin2023}
E.~Gardella, S.~Geffen, J.~Kranz, and P.~Naryshkin.
\newblock Classifiability of crossed products by nonamenable groups.
\newblock {\em J. Reine Angew. Math.}, 797:285--312, 2023.

\bibitem[GKPT18]{GKPT}
T.~Giordano, D.~Kerr, N.~C. Phillips, and A.~S. Toms.
\newblock {\em Crossed Products of {$C^*$}-Algebras, Topological Dynamics, and
  Classification}.
\newblock Advanced Courses in Mathematics -- CRM Barcelona. Birkh{\"a}user,
  Cham, 2018.

\bibitem[GL26]{glasnerliu}
E.~Glasner and C.~Liu.
\newblock Comparison and almost finiteness for actions of amenable groups,
  2026.
\newblock Preprint, arXiv:2609.23702.

\bibitem[Gla03]{Glasner_ergodicity_joining}
E.~Glasner.
\newblock {\em Ergodic theory via joinings}, volume 101 of {\em Mathematical
  Surveys and Monographs}.
\newblock American Mathematical Society, Providence, RI, 2003.

\bibitem[GLN20a]{GLN1}
G.~Gong, H.~Lin, and Z.~Niu.
\newblock A classification of finite simple amenable {$\mathcal{Z}$}-stable
  {$C^\ast$}-algebras, {I}: {$C^\ast$}-algebras with generalized tracial rank
  one.
\newblock {\em C. R. Math. Acad. Sci. Soc. R. Can.}, 42(3):63--450, 2020.

\bibitem[GLN20b]{GLN2}
G.~Gong, H.~Lin, and Z.~Niu.
\newblock A classification of finite simple amenable {$\mathcal{Z}$}-stable
  {$C^\ast$}-algebras, {II}: {${ C}^\ast$}-algebras with rational generalized
  tracial rank one.
\newblock {\em C. R. Math. Acad. Sci. Soc. R. Can.}, 42(4):451--539, 2020.

\bibitem[GPS95]{GPS}
T.~Giordano, I.~F. Putnam, and C.~F. Skau.
\newblock Topological orbit equivalence and {$C^*$}-crossed products.
\newblock {\em J. Reine Angew. Math.}, 469:51--111, 1995.

\bibitem[GW79]{glasner1979construction}
S.~Glasner and B.~Weiss.
\newblock On the construction of minimal skew products.
\newblock {\em Israel J. Math.}, 34(4):321--336, 1979.

\bibitem[Hal54]{hall}
P.~Hall.
\newblock Finiteness conditions for soluble groups.
\newblock {\em Proc. London Math. Soc. (3)}, 4:419--436, 1954.

\bibitem[Has16]{hasegawa16}
Kei Hasegawa.
\newblock {$KK$}-equivalence for amalgamated free product {$ C^*$}-algebras.
\newblock {\em Int. Math. Res. Not. IMRN}, (24):7619--7636, 2016.

\bibitem[Hat02]{Hatcher}
A.~Hatcher.
\newblock {\em Algebraic topology}.
\newblock Cambridge University Press, Cambridge, 2002.

\bibitem[Hig40]{Higman40}
G.~Higman.
\newblock The units of group-rings.
\newblock {\em Proc. London Math. Soc.}, 46(1):231--248, 1940.

\bibitem[Hig51]{Higman1951}
G.~Higman.
\newblock A finitely generated infinite simple group.
\newblock {\em J. London Math. Soc.}, 26(1):61--64, 1951.

\bibitem[HK01]{HK-BC}
N.~Higson and G.~Kasparov.
\newblock {$E$}-theory and {$KK$}-theory for groups which act properly and
  isometrically on {H}ilbert space.
\newblock {\em Invent. Math.}, 144(1):23--74, 2001.

\bibitem[HP26]{hp}
I.~Hirshberg and N.~C. Phillips.
\newblock An isomorphism theorem for infinite reduced free products, 2026.
\newblock Preprint, arXiv:2602.10220.

\bibitem[HW26]{huawhite}
S.~Hua and S.~White.
\newblock Uniqueness for embeddings of nuclear {$C^*$}-algebras into type
  {II}$_{1}$ factors, 2026.
\newblock Preprint, arXiv:2601.08779.

\bibitem[IPV13]{IPV13}
A.~Ioana, S.~Popa, and S.~Vaes.
\newblock A class of superrigid group von {N}eumann algebras.
\newblock {\em Ann. of Math. (2)}, 178(1):231--286, 2013.

\bibitem[ISTW26]{OWR2025Cstar}
A.~Ioana, D.~Shlyakhtenko, A.~Thom, and W.~Winter.
\newblock ${C}^*$-algebras.
\newblock {\em Oberwolfach Reports}, 22(3):1887--1948, 2026.

\bibitem[JS99]{Jiang-Su}
X.~Jiang and H.~Su.
\newblock On a simple unital projectionless {$C^*$}-algebra.
\newblock {\em Amer. J. Math.}, 121(2):359--413, 1999.

\bibitem[Kai04]{Kaimanovich2004}
V.~A. Kaimanovich.
\newblock Boundary amenability of hyperbolic spaces.
\newblock In {\em Discrete Geometric Analysis}, volume 347 of {\em Contemporary
  Mathematics}, pages 83--111. American Mathematical Society, Providence, RI,
  2004.

\bibitem[Kal69]{kallman}
R.~R. Kallman.
\newblock A generalization of free action.
\newblock {\em Duke Math. J.}, 36:781--789, 1969.

\bibitem[KB02]{KapovichBenakli2002}
I.~Kapovich and N.~Benakli.
\newblock Boundaries of hyperbolic groups.
\newblock In {\em Combinatorial and geometric group theory (New York,
  2000/Hoboken, NJ, 2001)}, volume 296 of {\em Contemporary Mathematics}, pages
  39--93. American Mathematical Society, Providence, RI, 2002.

\bibitem[Kel31]{Keller1931}
O.-H. Keller.
\newblock {Die Homoiomorphie der kompakten konvexen Mengen im Hilbertschen
  Raum}.
\newblock {\em Math. Ann.}, 105:748--758, 1931.

\bibitem[KES25]{SriSchafhauser}
S.~Kunnawalkam~Elayavalli and C.~Schafhauser.
\newblock Negative resolution to the {$C^*$}-algebraic {Tarski} problem, 2025.
\newblock Preprint, arXiv:2503.10505.

\bibitem[Kir95]{Kirchberg_classi}
E.~Kirchberg.
\newblock Exact {${ C}^*$}-algebras, tensor products, and the classification of
  purely infinite algebras.
\newblock In {\em Proceedings of the {I}nternational {C}ongress of
  {M}athematicians, {V}ol.\ 1, 2 ({Z}\"urich, 1994)}, pages 943--954.
  Birkh\"auser, Basel, 1995.

\bibitem[Kis81]{kishimoto}
A.~Kishimoto.
\newblock Outer automorphisms and reduced crossed products of simple
  {$C\sp{\ast} $}-algebras.
\newblock {\em Comm. Math. Phys.}, 81(3):429--435, 1981.

\bibitem[KK17]{KalantarKennedy2017}
M.~Kalantar and M.~Kennedy.
\newblock Boundaries of reduced {$C^*$}-algebras of discrete groups.
\newblock {\em J. Reine Angew. Math.}, 727:247--267, 2017.

\bibitem[KL26]{David_Li_2024}
D.~Kerr and H.~Li.
\newblock The small boundary property in products.
\newblock {\em Proc. Amer. Math. Soc.}, 154(10):4219--4231, 2026.

\bibitem[KN25]{KerrNaryshkin}
D.~Kerr and P.~Naryshkin.
\newblock Elementary amenability and almost finiteness.
\newblock {\em Compositio Math.}, 161(12):3321--3337, 2025.

\bibitem[KP00]{KP_classi}
E.~Kirchberg and N.~C. Phillips.
\newblock Embedding of exact {$C^*$}-algebras in the {C}untz algebra
  {$\mathcal{O}_2$}.
\newblock {\em J. Reine Angew. Math.}, 525:17--53, 2000.

\bibitem[KR14]{KR_SI}
E.~Kirchberg and M.~R{\o }rdam.
\newblock Central sequence {$C^*$}-algebras and tensorial absorption of the
  {J}iang-{S}u algebra.
\newblock {\em J. Reine Angew. Math.}, 695:175--214, 2014.

\bibitem[KT76]{KanThurston1976}
D.~M. Kan and W.~P. Thurston.
\newblock Every connected space has the homology of a {$K(\pi,1)$}.
\newblock {\em Topology}, 15(3):253--258, 1976.

\bibitem[Laf12]{lafforgue12}
V.~Lafforgue.
\newblock La conjecture de {B}aum-{C}onnes \`a{} coefficients pour les groupes
  hyperboliques.
\newblock {\em J. Noncommut. Geom.}, 6(1):1--197, 2012.

\bibitem[Lea13]{Leary2013}
I.~J. Leary.
\newblock A metric {Kan--Thurston} theorem.
\newblock {\em J. Topol.}, 6(1):251--284, 2013.

\bibitem[Li20]{Li_Cartan}
X.~Li.
\newblock Every classifiable simple {$\rm C^*$}-algebra has a {C}artan
  subalgebra.
\newblock {\em Invent. Math.}, 219(2):653--699, 2020.

\bibitem[LS77]{lyndon-schupp}
R.~C. Lyndon and P.~E. Schupp.
\newblock {\em Combinatorial Group Theory}, volume~89 of {\em Ergebnisse der
  Mathematik und ihrer Grenzgebiete}.
\newblock Springer-Verlag, Berlin--New York, 1977.

\bibitem[MN06]{meyernest}
Ralf Meyer and Ryszard Nest.
\newblock The {B}aum-{C}onnes conjecture via localisation of categories.
\newblock {\em Topology}, 45(2):209--259, 2006.

\bibitem[Mun00]{Munkres2000}
J.~R. Munkres.
\newblock {\em Topology}.
\newblock Prentice Hall, Upper Saddle River, NJ, second edition, 2000.

\bibitem[MvN36]{mvn36}
F.~J. Murray and J.~von Neumann.
\newblock On rings of operators.
\newblock {\em Ann. of Math. (2)}, 37(1):116--229, 1936.

\bibitem[MW24]{mawu}
X.~Ma and J.~Wu.
\newblock Almost elementary groupoid models for ${C}^*$-algebras, 2024.
\newblock Preprint, arXiv:2407.05251.

\bibitem[MY02]{MineyevYu2002}
I.~Mineyev and G.~Yu.
\newblock The {Baum--Connes} conjecture for hyperbolic groups.
\newblock {\em Invent. Math.}, 149:97--122, 2002.

\bibitem[Nea24]{RN_rr0}
R.~Neagu.
\newblock The admissible {KMS} bundles on classifiable {$C^*$}-algebras.
\newblock {\em Proc. Roy. Soc. Edinburgh Sect. A}, pages 1--34, 2024.

\bibitem[NR16]{ngrobert}
P.~W. Ng and L.~Robert.
\newblock Sums of commutators in pure {$C^*$}-algebras.
\newblock {\em M\"unster J. Math.}, 9(1):121--154, 2016.

\bibitem[Oml20]{omland}
T.~Omland.
\newblock Free nilpotent groups are {$ C^*$}-superrigid.
\newblock {\em Proc. Amer. Math. Soc.}, 148(1):283--287, 2020.

\bibitem[OO01]{oyono2001baum}
H.~Oyono-Oyono.
\newblock {Baum--Connes} conjecture and group actions on trees.
\newblock {\em K-Theory}, 24(2):115--134, 2001.

\bibitem[{Ope}26a]{openai-fgfp}
{OpenAI}.
\newblock An isomorphism of the free group factors.
\newblock OpenAI Math Release, September 2026.
\newblock Preprint, 23 September 2026.

\bibitem[{Ope}26b]{OpenAInonrigidity}
{OpenAI}.
\newblock Ten advances in mathematics and theoretical computer science.
\newblock OpenAI preprint, 2026.

\bibitem[{Ope}26c]{OpenAIrigidity}
{OpenAI}.
\newblock {The arithmetic category and stable recovery of lattice factors}.
\newblock OpenAI Math Release preprint, 2026.

\bibitem[Phi00]{Phillips2000}
N.~C. Phillips.
\newblock A classification theorem for nuclear purely infinite simple
  ${C}^*$-algebras.
\newblock {\em Doc. Math.}, 5:49--114, 2000.

\bibitem[Pop07]{popa07-rigidity}
S.~Popa.
\newblock Cocycle and orbit equivalence superrigidity for malleable actions of
  {$w$}-rigid groups.
\newblock {\em Invent. Math.}, 170(2):243--295, 2007.

\bibitem[Ren08]{Renault_Cartan}
J.~Renault.
\newblock Cartan subalgebras in {$C^*$}-algebras.
\newblock {\em Irish Math. Soc. Bull.}, (61):29--63, 2008.

\bibitem[Rob12]{robertNCCW}
L.~Robert.
\newblock Classification of inductive limits of 1-dimensional {NCCW} complexes.
\newblock {\em Adv. Math.}, 231(5):2802--2836, 2012.

\bibitem[Rob25]{robert}
L.~Robert.
\newblock Selfless {$C^*$}-algebras.
\newblock {\em Adv. Math.}, 478:Paper No. 110409, 28, 2025.

\bibitem[RS02]{RordamStormer}
M.~R{\o}rdam and E.~St{\o}rmer.
\newblock {\em Classification of Nuclear {$C^*$}-Algebras. Entropy in Operator
  Algebras}, volume 126 of {\em Encyclopaedia of Mathematical Sciences}.
\newblock Springer-Verlag, Berlin, 2002.

\bibitem[Sat19]{sato}
Y.~Sato.
\newblock Actions of amenable groups and crossed products of
  {$\mathcal{Z}$}-absorbing {$ C^*$}-algebras.
\newblock In {\em Operator algebras and mathematical physics}, volume~80 of
  {\em Adv. Stud. Pure Math.}, pages 189--210. Math. Soc. Japan, Tokyo, 2019.

\bibitem[Sch74]{scheinberg}
S.~Scheinberg.
\newblock Homeomorphism and isomorphism of abelian groups.
\newblock {\em Canadian J. Math.}, 26:1515--1519, 1974.

\bibitem[Sch20a]{new-TWW}
C.~Schafhauser.
\newblock A new proof of the {T}ikuisis-{W}hite-{W}inter theorem.
\newblock {\em J. Reine Angew. Math.}, 759:291--304, 2020.

\bibitem[Sch20b]{Schafhauser.Annals}
C.~Schafhauser.
\newblock Subalgebras of simple {AF}-algebras.
\newblock {\em Ann. of Math. (2)}, 192(2):309--352, 2020.

\bibitem[Shl26]{shlyakhtenko-semicircular}
D.~Shlyakhtenko.
\newblock On the isomorphism of semicircular algebras, 2026.
\newblock Preprint, arXiv:2610.19590.

\bibitem[Spi07a]{Spielberg_graph_model}
J.~Spielberg.
\newblock Graph-based models for {K}irchberg algebras.
\newblock {\em J. Operator Theory}, 57(2):347--374, 2007.

\bibitem[Spi07b]{Spielberg_prime_order}
J.~Spielberg.
\newblock Non-cyclotomic presentations of modules and prime-order automorphisms
  of {K}irchberg algebras.
\newblock {\em J. Reine Angew. Math.}, 613:211--230, 2007.

\bibitem[STW26]{99_open_nuclear}
C.~Schafhauser, A.~Tikuisis, and S.~White.
\newblock Nuclear ${C}^*$-algebras: 99 problems, 2026.
\newblock To appear in Special Issue of Munster J. Math. in memory of Eberhard
  Kirchberg.

\bibitem[Suz17]{suzuki2017construction}
Y.~Suzuki.
\newblock Construction of minimal skew products of amenable minimal dynamical
  systems.
\newblock {\em Groups Geom. Dyn.}, 11(1):75--94, 2017.

\bibitem[Sza19]{Gabor_actions_SSA}
G.~Szab\'o.
\newblock Actions of certain torsion-free elementary amenable groups on
  strongly self-absorbing {$ C^*$}-algebras.
\newblock {\em Comm. Math. Phys.}, 371(1):267--284, 2019.

\bibitem[Sza26]{gabor}
G.~Szabó.
\newblock The uniqueness theorem for {K}asparov theory, 2026.
\newblock Preprint, arXiv:2601.23029.

\bibitem[Tho03]{thomsen03}
Klaus Thomsen.
\newblock On the {$KK$}-theory and the {$E$}-theory of amalgamated free
  products of {$C^\ast$}-algebras.
\newblock {\em J. Funct. Anal.}, 201(1):30--56, 2003.

\bibitem[Tu99]{tu}
J.-L. Tu.
\newblock La conjecture de {B}aum-{C}onnes pour les feuilletages moyennables.
\newblock {\em $K$-Theory}, 17(3):215--264, 1999.

\bibitem[TW07]{SSA}
A.~S. Toms and W.~Winter.
\newblock Strongly self-absorbing {$C^*$}-algebras.
\newblock {\em Trans. Amer. Math. Soc.}, 359(8):3999--4029, 2007.

\bibitem[TW09]{tomswinter}
A.~S. Toms and W.~Winter.
\newblock Minimal dynamics and the classification of {$C^*$}-algebras.
\newblock {\em Proc. Natl. Acad. Sci. USA}, 106(40):16942--16943, 2009.

\bibitem[TWW17]{TWW}
A.~Tikuisis, S.~White, and W.~Winter.
\newblock Quasidiagonality of nuclear {$C^\ast$}-algebras.
\newblock {\em Ann. of Math. (2)}, 185(1):229--284, 2017.

\bibitem[Urs21]{ursu}
D.~Ursu.
\newblock Characterizing traces on crossed products of noncommutative
  {$C^*$}-algebras.
\newblock {\em Adv. Math.}, 391:Paper No. 107955, 29, 2021.

\bibitem[Val02]{valette}
A.~Valette.
\newblock {\em Introduction to the {B}aum-{C}onnes conjecture}.
\newblock Lectures in Mathematics ETH Z\"urich. Birkh\"auser Verlag, Basel,
  2002.
\newblock From notes taken by Indira Chatterji, With an appendix by Guido
  Mislin.

\bibitem[Var85]{Varadarajan1985}
K.~Varadarajan.
\newblock Pseudo-mitotic groups.
\newblock {\em J. Pure Appl. Algebra}, 37(2):205--213, 1985.

\bibitem[Wei94]{weibel}
C.~A. Weibel.
\newblock {\em An introduction to homological algebra}, volume~38 of {\em
  Cambridge Studies in Advanced Mathematics}.
\newblock Cambridge University Press, Cambridge, 1994.

\bibitem[Whi23]{stuart}
S.~White.
\newblock Abstract classification theorems for amenable {$C^\ast$}-algebras.
\newblock In {\em I{CM}---{I}nternational {C}ongress of {M}athematicians.
  {V}ol. 4. {S}ections 5--8}, pages 3314--3338. EMS Press, Berlin, 2023.

\bibitem[Wil07]{williams}
D.~P. Williams.
\newblock {\em Crossed products of {$C{^\ast}$}-algebras}, volume 134 of {\em
  Mathematical Surveys and Monographs}.
\newblock American Mathematical Society, Providence, RI, 2007.

\bibitem[Wu25]{Wu}
V.~Wu.
\newblock {$C^*$}-algebras associated to directed graphs of groups, and models
  of {K}irchberg algebras.
\newblock {\em J. Funct. Anal.}, 288(3):Paper No. 110740, 39, 2025.

\bibitem[Zho26]{zhou}
S.~Zhou.
\newblock {ICC} property {(T)} groups without {$W^*$}-superrigidity, 2026.

\bibitem[ZM68]{zeller-meier}
G.~Zeller-Meier.
\newblock Produits crois\'es d'une {$C\sp{\ast} $}-alg\`ebre par un groupe
  d'automorphismes.
\newblock {\em J. Math. Pures Appl. (9)}, 47:101--239, 1968.

\end{thebibliography}
\end{document}